\documentclass[a4paper]{article}

\usepackage{arxiv}

\usepackage{bookmark}
\usepackage{booktabs}
\usepackage{hyperref}
\usepackage{graphicx}
\usepackage{subfigure}
\usepackage{amsmath,amsthm,amssymb}
\usepackage{amsfonts}
\usepackage{enumerate}
\usepackage[nosort]{cite}
\usepackage{url}
\usepackage{xcolor}
\usepackage{tikz}
\usepackage{siunitx,cleveref}
\usepackage{mathtools}

\usetikzlibrary{graphs, quotes, positioning, shapes.geometric, arrows.meta}

\numberwithin{equation}{section}

\newtheorem{remark}{Remark}

\crefname{section}{Section}{Sections}
\crefname{figure}{Figure}{Figures}
\crefname{table}{Table}{Tables}

\newcommand{\R}{\mathbb{R}}

\title{An Acoustic Inversion-Based Flow Measurement Model in 3D Hydrodynamic Systems}

\renewcommand{\shorttitle}{Feasibility study of flow measurement method}

\author{Jiwei Li$^{1}$, Lingyun Qiu$^{2,3}$, Zhongjing Wang$^{4,5}$, Hui Yu$^{6}$\\
$^1$ Department of Mathematical Sciences, Tsinghua University, Beijing 100084, China\\
(\href{mailto:li-jw20@mails.tsinghua.edu.cn}{li-jw20@mails.tsinghua.edu.cn})\\
$^2$ Yau Mathematical Sciences Center, Tsinghua University, Beijing 100084, China\\
$^3$ Yanqi Lake Beijing Institute of Mathematical Sciences and Applications, Beijing 101408, China\\
(\href{mailto:lyqiu@tsinghua.edu.cn}{lyqiu@tsinghua.edu.cn})\\
$^4$ Department of Hydraulic Engineering, Tsinghua University, Beijing 100084, China\\
$^5$ Breeding Base for State Key Laboratory of Land Degradation and Ecological Restoration in Northwest China, \\
Ningxia University, Yinchuan 750021, Ningxia, China \\
(\href{mailto:zj.wang@tsinghua.edu.cn}{zj.wang@tsinghua.edu.cn})\\
$^6$ School of Mathematics and Computational Science, Xiangtan University, Xiangtan 411105, Hunan, China\\
(\href{mailto:huiyu@xtu.edu.cn}{huiyu@xtu.edu.cn})
}

\begin{document}
\maketitle

\begin{abstract}
	This study extends the flow measurement method initially proposed in \cite{liFlowMeasurementInverse2023} to three-dimensional scenarios, addressing the growing need for accurate and efficient non-contact flow measurement techniques in complex hydrodynamic environments. Compared to conventional Acoustic Doppler Current Profilers (ADCPs) and remote sensing-based flow monitoring, the proposed method enables high-resolution, continuous water velocity measurement, making it well-suited for hazardous environments such as floods, strong currents, and sediment-laden rivers. Building upon the original approach, we develop an enhanced model that incorporates multiple emission directions and flexible configurations of receivers. These advancements improve the adaptability and accuracy of the method when applied to three-dimensional flow fields. To evaluate its feasibility, we conducted extensive numerical simulations designed to mimic real-world hydrodynamic conditions. The results demonstrate that the proposed method effectively handles diverse and complex flow field configurations, highlighting its potential for practical applications in water resource management and hydraulic engineering.
\end{abstract}



\section{Introduction}\label{sec_intro}

Fluid flow is a central topic in fluid mechanics serving as a cornerstone of physics and engineering. It manifests across a broad spectrum of scales and environments: large-scale flows such as open channel flow and river systems, as well as small-scale flows like blood circulation and water distribution networks. These flows play a crucial role in both natural processes and engineered systems. However, not all fluid flows are inherently beneficial. For instance, uncontrolled river flows can result in devastating floods, damaging infrastructure and crops, while irregular blood flow in vessels can lead to severe health risks and potentially life-threatening outcomes\cite{kheyfetsAbstract12086Irregular2014}.

The ability to monitor, control, and utilize various fluid flows is of importance to humanity. Effective management of flows not only helps to mitigate disasters, such as floods and blockages, but also enhances the efficiency of flow utilization in diverse applications. However, the invisible and intangible nature of fluid flow presents significant challenges for its measurement and control. As a result, reliable flow monitoring techniques\cite{tauroNovelPermanentGaugecam2016,sahooMultimissionVirtualMonitoring2024,burauInnovationMonitoringUS2016,walkerEphemeralSandRiver2019,luMethodMonitoringEnvironmental2022,johansenMonitoringCoastalWater2022,bjerklieSatelliteRemoteSensing2018}, ranging from flow gauging stations and hydrological prediction to satellite remote sensing, form the essential foundation for flow management. Among these, accurate flow measurement stands out as a critical component, enabling precise control\cite{zengPreciseMeasurementControl2020} and informed decision-making\cite{padikkalInformedDecisionMaking2013} in both natural and engineered systems.

Flow measurement is a critical process across diverse fields, including physics, engineering, and environmental sciences. Key flow characteristics such as velocity, pressure, and temperature are essential for understanding fluid dynamics, with flow velocity often being the most significant parameter due to its direct relationship with the rate of fluid movement. However, accurately measuring flow in complex environments, such as rivers with varying cross-sections, lakes with irregular boundaries, and oceans with dynamic currents, remains a significant challenge. Traditional methods often struggle with these complexities, especially in confined or turbulent conditions. Numerous techniques exist for measuring flow velocity, which can be broadly categorized into two groups: intrusive and non-intrusive methods.

Intrusive methods, such as pressure-difference-based techniques (e.g., Pitot tubes, Venturi meters, and orifice plates, see \cite{liptakInstrumentEngineersHandbook2003c,blincoTurbulenceCharacteristicsFree1971,klopfensteinjrAirVelocityFlow1998} for more details), require the insertion of devices into the flow path, potentially altering the flow conditions. In contrast, non-intrusive methods, including laser Doppler velocimetry\cite{georgeLaserDopplerVelocimeterIts1973,meierImagingLaserDoppler2012}, acoustic Doppler velocimetry\cite{songTurbulenceMeasurementNonuniform2001,garciaTurbulenceMeasurementsAcoustic2005}, electromagnetic meters\cite{shercliffTheoryElectromagneticFlowMeasurement1962,watralElectromagneticFlowMeters2015}, thermal anemometry\cite{orluThermalAnemometry2017,fossThermalTransientAnemometer2004}, Coriolis meters\cite{anklinCoriolisMassFlowmeters2006,wangCoriolisFlowmetersReview2014}, and particle image velocimetry\cite{wereleyRecentAdvancesMicroParticle2010,westerweelParticleImageVelocimetry2013,willertDigitalParticleImage1991,raffelParticleImageVelocimetry2018}, measure flow characteristics without physically disturbing the flow. Due to their minimal impact on the system, non-intrusive methods have become increasingly popular in practical applications. This study focuses on acoustic wave-based techniques, particularly ADCPs, as a promising approach for accurate and efficient flow velocity measurement in three-dimensional scenarios.

Li et al.\ introduced a novel flow measurement model inspired by ADCPs in \cite{liFlowMeasurementInverse2023}. This approach leverages the inverse source problem \cite{isakovInverseSourceProblems1990} of wave equations to detect the distribution of particles within the flow at any given time. By analyzing a sequence of particle distributions, the flow field can then be reconstructed using various computational techniques. The proposed model demonstrated strong performance in both theoretical analysis and numerical simulations. However, its evaluation was limited to simplified, two-dimensional environments in \cite{liFlowMeasurementInverse2023}. In this study, we extend the investigation to explore the feasibility of this model in three-dimensional settings that are more complex. Through comprehensive numerical simulations, we evaluate and validate performances of the model in realistic hydrodynamic scenarios, highlighting its potential for broader applications.

Sidelobe reflections \cite{lentzNoteDepthSidelobe2022} are well known to significantly impact the accuracy of flow velocity measurements obtained using ADCPs, often introducing substantial errors in various scenarios. Although the acoustic sidelobes are considerably weaker than the main beams \cite{appellAcousticDopplerCurrent1991,gordonAcousticDopplerCurrent1996}, surface reflections can be up to 100 times stronger than those from particles suspended in the water. As a result, sidelobe reflections from the surface can rival the strength of the main beam returns from the water column. This means that acoustic returns from the main beams are only reliable for calculating currents before the vertical sidelobe reflections from the sea surface arrive \cite{appellAcousticDopplerCurrent1991,gordonAcousticDopplerCurrent1996,lentzNoteDepthSidelobe2022}. These limitations significantly restrict the applicability of ADCPs in scenarios where strong surface reflections are unavoidable. The method proposed in \cite{liFlowMeasurementInverse2023} has shown potential for measuring flow velocity even in the presence of strong sidelobe reflections. However, that study did not address the impact of boundary reflections, which often arise in confined or complex environments and pose additional challenges for accurate measurement. In this paper, we extend the previous model to explicitly consider the effects of boundary reflections, incorporating these factors into our numerical simulations. This enhancement allows us to evaluate the performance of the method under more realistic conditions. Moreover, the proposed approach is not limited to simple configurations such as canals. It can be applied to a variety of complex flow environments, including rivers, lakes, and oceans, where the presence of irregular and dynamic boundary conditions adds further complexity. By addressing these challenges, our method demonstrates greater adaptability and reliability in diverse hydrodynamic systems.

The paper is structured as follows. \Cref{sec_theory} provides an introduction to the theory and methods of the flow measurement approach based on the inverse source problem. In \cref{sec_numerical}, we present a series of numerical results demonstrating the application of this method in various flow environments. Finally, \cref{sec_conclusion} concludes the paper, summarizing the key findings and potential future directions.

\section{Theory and methods}\label{sec_theory}
In this section, we extend the two-dimensional flow measurement model based on the inverse source problem to three-dimensional space $\R^3$, while incorporating the effects of sidelobe reflections on the accuracy and stability of the velocity estimation process. For simplicity, we refer to the proposed flow measurement model as the Acoustic Inversion-based Flow Measurement (AIFM) model. The discussion begins with a review of the fundamental measurement workflow of the AIFM model, followed by its three-dimensional extension. Each step in the workflow is then examined in detail, highlighting the specific challenges and implementations required in three-dimensional scenarios.

The AIFM model's overall measurement workflow is summarized as follows:
\begin{enumerate}
	\item Acoustic transducers emit wave signals into the flow field. These waves propagate through the fluid and interact with moving particles, generating scattered signals.
	\item Receivers collect the scattered acoustic signals and use the recorded pressure data as input for further processing. This processing involves solving the inverse source problem and calculating the velocity field.
\end{enumerate}

As noted in \cref{sec_intro}, sidelobe reflections significantly impact the precision of velocity field estimation. Although the mathematical foundation and numerical algorithms of the AIFM model in two-dimensional space $\R^2$, have been thoroughly validated, real-world flow measurement is inherently three-dimensional and constrained within finite domains. Furthermore, diverse boundary conditions in such domains critically affect the model's stability and feasibility.

To address these challenges, this work extends the AIFM model to three-dimensional space $\R^3$. More specifically, this involves rigorous theoretical analysis to adapt the inverse source problem from an unbounded domain to bounded regions, accounting for various boundary conditions and their influence on the uniqueness and stability of the solution. Additionally, specific adjustments to the measurement workflow are necessary to accommodate three-dimensional propagation effects and ensure the robustness of the velocity estimation under realistic conditions. By explicitly incorporating sidelobe reflections, the extended model bridges the gap between theoretical results and practical applications, enabling reliable velocity field measurements even in complex environments.

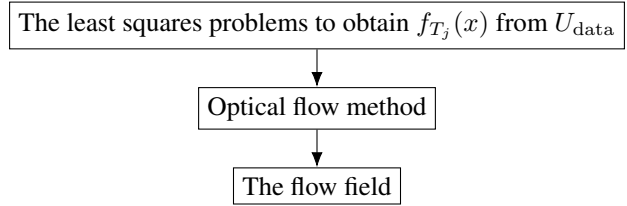
\begin{figure}
	\centering
	\begin{tikzpicture}[node distance=5mm
		]
		\node[draw,align=center]	(solve)		{The least squares problems to obtain $f_{T_j}(x)$ from $U_\mathrm{data}$};
		\node[draw,below=of solve,align=center]		(of)	{ Optical flow method};
		\node[draw,below= of of]		(flow)	{The flow field};
		\draw[-{Latex[length=2mm]}] (solve) -- (of);
		\draw[-{Latex[length=2mm]}] (of) -- (flow);
	\end{tikzpicture}
	\caption{The process diagram of inversion for the velocity of particles.}
	\label{fig_algorithm}
\end{figure}

The AIFM model's algorithmic processing workflow is illustrated in \cref{fig_algorithm} and comprises two key steps: particle detection and velocity field computation. In the first step, the distribution of moving particles in the fluid is reconstructed at each time point using the data received by boundary receivers. In the second step, the velocity field of the flow is calculated using optical flow techniques applied to the reconstructed particle distributions. These techniques include methods such as auto-correlation and cross-correlation.It is important to note that the velocity field obtained in the second step reflects the motion of particles within the fluid rather than the fluid itself. Consequently, the AIFM model assumes that the velocities of the particles are consistent with the velocity field of the fluid. This assumption is crucial to the successful application of the method in practical scenarios.

In the subsequent sections, we will present a detailed discussion of the three-dimensional extension of the AIFM model. Specifically, we will explore how each step in the measurement workflow is adapted for three-dimensional implementation, addressing the unique challenges posed by the inclusion of sidelobe reflections and the complexities of bounded domains. These considerations are essential for ensuring the accuracy and robustness of the model when applied to real-world fluid environments.

\subsection{Particle detection based on the AIFM method}\label{sec_isp}
In this section, we introduce the theoretical foundation of the AIFM method for particle detection in three-dimensional flow environments. The method leverages the inverse source problem of wave equations to reconstruct the spatial distribution of particles within the fluid. By analyzing the scattered acoustic signals collected by boundary receivers, the AIFM method infers the locations of moving particles at each time point, enabling accurate and efficient flow measurement.

We consider a three-dimensional strip domain $\Omega=[0,L_1]\times \R \times [0,L_3]\subset\R^3$, which models an open channel environment. Here, $L_1$ and $L_3$ denote the width and depth of the open channel, respectively. Let $x=(x_1,x_2,x_3)$, where $x_1$ represents the width direction of the open channel, $x_2$ represents the length direction of the open channel, and $x_3$ represents the vertical distance from position $x$ to the riverbed within the open channel. Let $\Gamma=\partial\Omega$ represent the boundary of $\Omega$, which is divided into two sub-boundaries: $\Gamma_1=\Gamma\cap\{x_3=L_3\}$, corresponding to the flow surface, and $\Gamma_2=\Gamma\setminus\Gamma_1$, representing the bed and sidewalls of the open channel.

The classification of $\Gamma$ into $\Gamma_1$ and $\Gamma_2$ accounts for the distinct acoustic reflection effects at liquid-gas and liquid-solid interfaces. At $\Gamma_1$, the interface between the liquid and gas, the acoustic pressure tends to zero when waves propagate from the liquid to the gas. This behavior is approximated using a Dirichlet boundary condition. Conversely, at $\Gamma_2$, representing the liquid-solid interface, the acoustic pressure remains constant along the normal direction of the interface. This property is modeled using a Neumann boundary condition.

This boundary condition treatment enables an accurate representation of wave propagation and reflection effects in open channel environments. It allows the model to incorporate sidelobe reflections, thereby enhancing the applicability and reliability of the AIFM method for measuring flow velocity in realistic scenarios. The three-dimensional extension of the AIFM model provides theoretical support for flow measurement in open channels and lays the groundwork for tackling more complex fluid systems.

We denote $f_t(x)$ as the distribution of moving particles within the fluid medium at time $t$, $U(x,t)$ as the acoustic pressure distribution of the wave signal at time $t$, and $\lambda(x,t)$ as the known emitted wave field. Considering the reflection effects on the two parts of the boundary, $\Gamma_1$ and $\Gamma_2$, as described above, the three-dimensional extension of the particle detection model in the AIFM framework is expressed as follows:
\begin{subequations}\label{eq_wave_system}
	\begin{align}
		 \frac{1}{c^2(x)} \partial_t^2 U(x,t)- \Delta U(x,t) = \lambda(x,t)f_t(x),& \quad (x,t) \in \Omega \times \R_+, \label{eq_wave_eqn} \\
		 U(x,t) =  0,& \quad (x,t) \in \Gamma_1 \times \R_+,\label{eq_wave_bc1} \\
		 \partial_\nu U(x,t) =  0,& \quad (x,t) \in \Gamma_2 \times \R_+,\label{eq_wave_bc2} \\
		 U(x,0) = 0,& \quad x \in \Omega, \label{eq_wave_ic1} \\
		 \partial_t U(x,0) = 0,& \quad x \in \Omega,\label{eq_wave_ic2}
	\end{align}
\end{subequations}
where $\R_+=(0,\infty)$ is the time domain, $\nu$ is the outward unit normal vector to the boundary $\Gamma_2$, and $c(x)$ is the sound speed of the flow at point $x$, typically around \SI[separate-uncertainty=true]{1500 \pm 10}{\meter\per\second}.

In traditional acoustic-based flow measurement techniques, signals reflected from boundaries are typically treated as noise, which can significantly degrade the accuracy of velocity field measurements\cite{lentzNoteDepthSidelobe2022}. In contrast, our approach explicitly incorporates boundary reflection effects into the flow measurement model. This allows for a precise description of acoustic wave propagation and reflection within open-channel environments. By doing so, we transform what was previously regarded as data contamination into an integral part of the model.

This innovative strategy not only improves the applicability and reliability of the AIFM model but also turns the influence of boundary reflections into useful information that enhances data interpretation. By leveraging the propagation and reflection characteristics of acoustic waves, the AIFM model achieves a significant improvement in the accuracy and stability of flow velocity measurements. This method provides a novel perspective and technical support for addressing flow measurement challenges in complex fluid systems.

Next, we apply a uniform discretization to the time domain, such that $T_j=j\Delta t$, where $j=0,1,2,\dots$. Within each time interval $[T_j,T_{j+1}]$, the source term $f_t(x)$ is assumed to be constant with respect to $t$, i.e. $f_t(x)\approx f_{T_j}(x)$, $t \in [T_j,T_{j+1}]$. Under this assumption, the original wave equation \eqref{eq_wave_system} can be approximated as follows:
\begin{subequations}\label{eq_approx_wave_system}
	\begin{align}
		 \frac{1}{c^2(x)} \partial_t^2 U(x,t)- \Delta U(x,t) = \lambda(x,t)f_{T_j}(x),& \quad (x,t) \in \Omega \times (T_j,T_{j+1}), \\
		 U(x,t) =  0,& \quad (x,t) \in \Gamma_1 \times (T_j,T_{j+1}),\label{eq_appro_wave_bc1}\\
		 \partial_\nu U(x,t) =  0,& \quad (x,t) \in \Gamma_2 \times (T_j,T_{j+1}),\label{eq_appro_wave_bc2}\\
		 U(x,T_j) = 0,& \quad x \in \Omega, \label{eq_appro_wave_ic1} \\
		 \partial_t U(x,T_j) = 0,& \quad x \in \Omega. \label{eq_appro_wave_ic2}
	\end{align}
\end{subequations} 
The validity of the approximated system \eqref{eq_approx_wave_system} can be demonstrated by combining reflection principle with the proof method in \cite{liFlowMeasurementInverse2023}. For conciseness, we omit the detailed proof here. 

To satisfy the initial conditions $U(x,T_j)=\partial_t U(x,T_j)=0$, it is necessary to select an appropriate measurement interval. This ensures that the acoustic wave field propagating within the previous measurement interval $[T_{j-1},T_j]$ decays to near zero within the measurement region at timedecays to near zero within the measurement region at time $T_j$. In practical computations, by selecting the measurement interval appropriately, this condition can be effectively met. Subsequently, the sequence of equations in \eqref{eq_approx_wave_system} can be summarized into a wave equation over a reference time interval $[0,T]$ as follows:
\begin{subequations}\label{eq_reduced_wave_system}
	\begin{align}
		 \frac{1}{c^2(x)} \partial_t^2 U(x,t)- \Delta U(x,t) = \lambda(x,t)f(x),& \quad (x,t) \in \Omega \times (0,T), \\
		 U(x,t) =  0,& \quad (x,t) \in \Gamma_1 \times (0,T),\label{eq_reduced_wave_bc1}\\
		 \partial_\nu U(x,t) =  0,& \quad (x,t) \in \Gamma_2 \times (0,T),\label{eq_reduced_wave_bc2}\\
		 U(x,0) = 0,& \quad x \in \Omega, \label{eq_reduced_wave_ic1} \\
		 \partial_t U(x,0) = 0,& \quad x \in \Omega. \label{eq_reduced_wave_ic2}
	\end{align}
\end{subequations}

Next, we will detail the location of the data collection area. Directly collecting data at the boundary $\Gamma$ of the measurement region is not suitable for several reasons: first, the data at $\Gamma$ is influenced by boundary conditions, which leads to a lack of received data on $\Gamma_1$; second, in practical scenarios, it is unrealistic to collect data exactly at the boundary $\Gamma$ due to physical constraints, such as the inability to place measurement devices directly at the water surface. Therefore, we need to collect data in a region near the boundary $\Gamma$ to reduce the impact of boundary effects. Specifically, we define the boundary location $\Gamma^\prime$ for data collection as part of the boundary of a subset $\Omega^\prime$ of $\Omega$, i.e. $\Gamma^\prime\subset\partial\Omega^\prime$, where $\Omega^\prime\subset\Omega$ satisfies that for any $t\in[0,T]$, the support of $f_t(\cdot)$ is contained within $\Omega^\prime$. In this way, we can better control the location and size of the data collection region, effectively reducing the interference of boundary effects on the data. In practical computations, the location and size of the data collection area can be selected appropriately to meet the specific application needs, further improving the accuracy and reliability of the measurement results.

The received data, denoted as $\{U(r_s,t)\}_{s=1}^N$, is collected at locations $r_s$ corresponding to $N$ receivers. These receivers are typically placed strategically on the boundary $\Gamma^\prime$ to capture the scattered signals generated by the interaction between the acoustic waves and the moving particles. This strategic placement maximizes the measurement range, enabling the capture of acoustic signals scattered from particles throughout the entire flow domain. By ensuring comprehensive coverage, this configuration enhances the reconstruction accuracy of the particle distribution.

Denote the acoustic signals collected at $r_s$ as $U_{\text{data},s}(t;\lambda)\coloneq U(r_s,t), s=1,2,\dots,N$, with the dependency on the emitted wave field $\lambda$ explicitly indicated. This dependency is crucial, as we will later discuss how varying the emitted wave field $\lambda$ can enhance the accuracy of the flow velocity measurement model. Then the primary goal of the AIFM method is to reconstruct the particle distribution $f(x)$ using the received data $U_{\text{data},s}$. Using this received data, we define the forward operator that relates the particle distribution $f(x)$ to the measured acoustic data $U_{\mathrm{data},s}(t;\lambda)$.
\begin{equation}\label{eq_forward_s-th}
	\mathcal{F}_s(f;\lambda)=U_{\text{data},s}(t;\lambda),\quad s=1,2,\dots,N.
\end{equation}
\begin{equation}\label{eq_forward}
	\mathcal{F}(f;\lambda)=\left(\mathcal{F}_1(f;\lambda),\mathcal{F}_2(f;\lambda),\dots,\mathcal{F}_N(f;\lambda)\right)
\end{equation}
The core of the AIFM method lies in reconstructing the particle distribution $f(x)$ by determining the inverse operator $\mathcal{G}^{-1}$ that relates the received data $U_{\mathrm{data}}(t;\lambda)$ to $f(x)$. Specifically, the inverse problem involves solving for $f(x)$ from the operator equation $\mathcal{G}(f; \lambda) = U_{\mathrm{data}}(t;\lambda)$. The uniqueness of the solution to this inverse problem can be established using an approach similar to \cite{liFlowMeasurementInverse2023}, combined with the reflection principle. This requires a tailored treatment for the strip domain $\Omega$, ensuring that the boundary conditions and geometry are appropriately incorporated into the analysis.

To further enhance the accuracy of the AIFM method, we propose the use of multiple emitted source waves. This approach draws inspiration from the multi-source inversion techniques commonly employed in geophysical imaging, which significantly increase the informational richness of the received data and improve the stability of the inversion process. The practicality of this method lies in its implementation, where varying parameters such as the position and orientation of the transmitter allow the emission of distinct wave fields.

Let the emitted wave fields be represented by $\lambda_m(x,t)$ for $m=1,2,\dots,M$. These multiple wave fields generate a set of forward operators $\mathcal{F}(f;\lambda_m)$, corresponding to the received data for each wave field. By combining these forward operators, the reconstruction of the particle distribution $f(x)$ becomes more robust and accurate. The optimization problem associated with the AIFM model is then redefined to leverage the aggregated information from multiple sources, expressed as:
\begin{equation}\label{eq_optimization_multi_source_waves}
	\min_f \mathcal{J}(f)\coloneq \frac{1}{2}\sum_{m=1}^M\sum_{s=1}^N \|\mathcal{F}_s(f;\lambda_m)-U_{\text{data},s}(\cdot;\lambda_m)\|^2_{L^2(0,T)}.
\end{equation}
where $\mathcal{J}(f)$ is the objective function to be minimized, and the optimization problem aims to find the particle distribution $f(x)$ that best fits the received data $U_{\text{data},s}(\cdot;\lambda_m)$ for each emitted wave field $\lambda_m$. 

To solve the optimization problem \eqref{eq_optimization_multi_source_waves}, various optimization techniques can be utilized. In this study, we adopt the conjugate gradient method, which has been validated for its efficiency and performance in \cite{liFlowMeasurementInverse2023}. This approach is particularly effective in handling high-dimensional optimization problems and ensures reliable convergence to a local minimum. Its proven performance makes it a robust choice for the complex reconstructions required in three-dimensional environments.

Once the particle distributions $f_{T_j}(x)$ are successfully reconstructed by solving the optimization problem \eqref{eq_optimization_multi_source_waves}, the first stage of the AIFM model is completed. The next stage focuses on computing the velocity field of the flow based on the reconstructed particle distribution. While the velocity field computation in the AIFM model shares conceptual similarities with the two-dimensional approach discussed in \cite{liFlowMeasurementInverse2023}, the complexities of a three-dimensional measurement domain necessitate more efficient and robust methodologies.

To address these challenges, we employ the optical flow method, which offers greater adaptability in three-dimensional spaces compared to traditional local velocity computation techniques. The optical flow method is well-suited for handling complex boundary conditions and non-uniform fluid motion, making it a powerful tool for accurate velocity field estimation in three-dimensional environments.

In the following section, we will provide a detailed discussion of the three-dimensional extension of the AIFM velocity computation model. We will focus on the implementation of the optical flow method in $\R^3$, analyze its performance, and validate its applicability in complex measurement domains through numerical experiments. The results will demonstrate the model's effectiveness in improving measurement precision and stability.

\subsection{Optical flow method}\label{sec_opticalFlow}
After the particle distributions $f_{T_j}(x)$ are obtained for each time step, the AIFM method employs the optical flow technique to compute the velocity field of the fluid. Given the computational challenges of three-dimensional data, we adopt the Farneb\"{a}ck method \cite{farnebackTwoframeMotionEstimation2003}, a robust optical flow algorithm based on the Taylor series expansion of image intensity functions. This method is renowned for its high accuracy and resilience to noise, making it well-suited for reconstructing velocity fields in complex flow scenarios.

To further improve computational efficiency, we leverage the Python implementation \texttt{farneback3d}, which harnesses the power of GPU-based parallel computing through CUDA. This approach dramatically accelerates the processing of large-scale 3D data, ensuring timely and precise flow measurements. By combining the accuracy of the Farneb\"{a}ck method with the computational speed of \texttt{farneback3d}, the AIFM model achieves a practical and efficient solution for three-dimensional velocity field estimation, paving the way for its application in real-world engineering scenarios.

\begin{remark}
	Note that the second stage of AIFM method is similar to the technique of particle image velocimetry (PIV, see \cite{scaranoTomographicPIVPrinciples2013, meinhartPIVMeasurementsMicrochannel1999, hainFundamentalsMultiframeParticle2007, elsingaTomographicParticleImage2006}), which also computes the entire velocity field by analyzing particle images at each time point. However, a key limitation of PIV is its inability to measure velocity components along the line of sight, which means additional cameras are needed to capture all components. In contrast, AIFM method can reconstruct the full velocity field with a single set of wave pressure data, without the need for multiple viewpoints. Furthermore, while PIV requires particles to be visible to lasers, the AIFM method only requires that particles be detectable by ultrasonic sound, a requirement easily met by most common materials. These advantages make the AIFM method a more versatile and practical choice for flow measurement in diverse environments.
\end{remark}

\section{Numerical results}\label{sec_numerical}
In this section, we present a series of numerical experiments to validate the feasibility and effectiveness of the AIFM method in three-dimensional flow environments. First, we examine the influence of various parameter selections on the accuracy and efficiency of particle detection using the AIFM framework. Next, the performances of the method in reconstructing flow velocity is evaluated across several representative flow fields. Finally, we demonstrate the practical applicability of the AIFM method by simulating flow fields commonly encountered in real-world scenarios. These experiments provide a comprehensive assessment of the robustness and adaptability of the method under diverse conditions.

\subsection{Particle detection}\label{subsec_particleDetection}
The performance of the AIFM method is influenced by various factors, including the frequency of the emitted waves, receiver configuration, number of receivers, particle density, emitted waveforms, and data noise. In this section, we conduct a comprehensive evaluation through numerical experiments to investigate how these factors affect the accuracy, robustness, and stability of the particle detection in AIFM method. Specifically, we examine the effects of number of receivers, particle density, and emitted waveforms on the reconstruction accuracy.

Previous studies, such as those in \cite{liFlowMeasurementInverse2023}, have demonstrated that the AIFM method exhibits strong adaptability and robustness to changes in wave frequency, receiver configurations, and noise levels. Building on these findings, this paper further explores the influence of additional parameters, such as the spatial density of particles, the diversity of emitted waveforms, and the total number of receivers, on the method’s detection accuracy. The results from these experiments provide valuable insights and technical guidance for optimizing the AIFM method, enabling better performance in complex measurement scenarios.

\subsubsection{Experimental setup}\label{subsubsec_numerical_setup}
In all numerical experiments conducted in this subsection, the width and depth of the open channel are set to $L_1=\SI{1}{\meter}$ and $L_3=\SI{1}{\meter}$, respectively. The sound speed $c(x)$ is assumed to be constant at \SI{1500}{\meter\per\second}. The spatial discretization step size is set to \SI{1}{\centi\meter}, while the time step size is chosen as \SI{2.3}{\micro\second} to ensure numerical stability and compliance with the Courant-Friedrichs-Lewy (CFL) condition. To achieve high accuracy in solving the wave equation, we utilize a high-order finite difference method. The implementation is automated using the Python library Devito\cite{louboutinDevitoV310Embedded2019,luporiniArchitecturePerformanceDevito2020}, which specializes in symbolic computation for stencil-based applications. Devito generates highly optimized finite difference schemes and leverages just-in-time compilation to produce efficient machine code, significantly enhancing computational efficiency.

For the emitted source wave $\lambda_m(x,t)$, we adopt the form of a traveling plane wave modeled as a Ricker wavelet\cite{ryanRickerOrmsbyKlauder1994}, which closely resembles real-world acoustic signals. The mathematical expression for the Ricker wavelet is as follows:
\begin{equation}\label{eq_ricker_wavelet}
	\lambda_m(x,t)=\left(1-2\pi^2 q_0^2 \left(t-\frac{p_m\cdot x}{c}-\frac{3}{\pi q_0}\right)^2\right)\cdot\exp\left(-\pi^2 q_0^2 \left(t-\frac{p_m\cdot x}{c}-\frac{3}{\pi q_0}\right)^2\right),
\end{equation}
where $q_0$ represents the central frequency of the Ricker wavelet, and $p_m$ denotes the direction of the $m$-th emitted wave. In the following experiments, the central frequency is fixed at \SI{20}{\kilo\hertz}. 

The choice of propagation direction $p_m$ significantly impacts the accuracy and stability of the AIFM model's particle detection process. A simple and intuitive method to choose these directions is to distribute them randomly on the unit sphere. This random selection ensures the isotropy of the wave propagation directions. However, the randomness of this approach leads to variability in the experimental results, reducing their reproducibility. Such variability can complicate the analysis and comparison of results, making it challenging to draw consistent conclusions.

\begin{figure}[ht]
	\centering
	\includegraphics[width=0.6\textwidth]{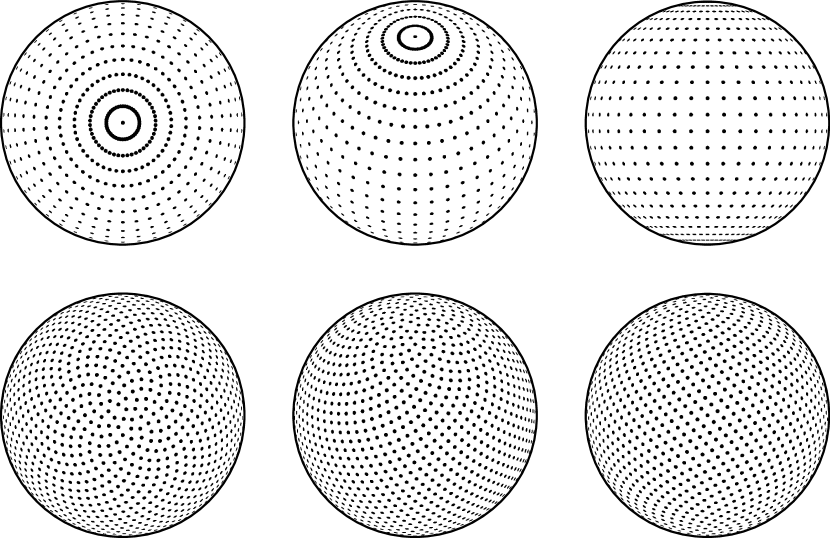}
	\caption{Comparison of the polar coordinate and Fibonacci lattice for selecting wave propagation directions\cite{gonzalezMeasurementAreasSphere2010}. The top row shows the polar coordinate distribution, while the bottom row displays the Fibonacci lattice distribution.}
	\label{fig_fibonacci_lattice}
\end{figure}

An alternative is to use polar coordinates for explicitly defining the directions. While this approach ensures repeatability, it suffers from uneven distribution, with directions clustering near the poles. To address this limitation, we adopt the Fibonacci lattice\cite{dixonSpiralPhyllotaxis1989,gonzalezMeasurementAreasSphere2010}. This method provide a quasi-uniform distribution of points on a sphere  and avoids the accumulation of emission directions at the poles, as illustrated in \cref{fig_fibonacci_lattice}. The Fibonacci lattice is not only computationally efficient but also ensures an even coverage of directions, thereby enhancing the diversity of the experimental data. The mathematical formulation for generating $M$ directions on the unit sphere is given by:
\begin{subequations}
	\begin{align}
		 p_{m,1}&=\sqrt{1-p_{m,3}^2}\cos(2\pi m \phi),\\
		 p_{m,2}&=\sqrt{1-p_{m,3}^2}\sin(2\pi m \phi),\\
		 p_{m,3}&=\frac{2m-1}{M}-1,
	\end{align}
\end{subequations}
where $\phi=(1+\sqrt{5})/2$ is the golden ratio, which ensures the near-uniform distribution of sampling points on the sphere. By leveraging this method, we can systematically select the propagation directions $p_m$, effectively addressing the drawbacks of previous approaches. The influences of these directions on the performance of the AIFM model will be further analyzed in subsequent numerical experiments to evaluate their effect on the detection model.

Spherical particles are randomly distributed within the domain, with diameters uniformly sampled from \SI{6}{\centi\meter} to \SI{10}{\centi\meter}, reflecting the typical sizes of particles found in practical flow scenarios, such as sediments or bubbles. The number of conjugate gradient iterations is fixed at 100, based on empirical observations during the solution of the inverse optimization problem \eqref{eq_optimization_multi_source_waves}. In practice, it has been observed that beyond 100 iterations, the reduction in the objective function $\mathcal{J}(f)$ becomes negligible, indicating that the optimization process has effectively converged. Thus, fixing the iteration count at 100 strikes a balance between ensuring solution accuracy and controlling computational costs. This setup provides a consistent baseline for evaluating the performance of the AIFM method across experiments involving variations in particle density, wave directions, and receiver configurations.

\subsubsection{Effect of the directions of emitted waves}\label{subsubsubsec_direction}
The choice of emitted wave directions significantly influences the performance of the AIFM method in particle detection. Specifically, the propagation directions $p_m$ of the Ricker wavelets \eqref{eq_ricker_wavelet} play a critical role in determining the stability and accuracy of the inversion results. This subsection examines the impact of different numbers of propagation directions on the performances of the AIFM model. By varying the number of wave propagation directions, we evaluate their effects on the model's detection stability and precision. Moreover, as the computational complexity of solving the inversion problem scales with the number of propagation directions, we investigate how this parameter affects computational efficiency. Through these experiments, we aim to provide a comprehensive understanding of the relationship between wave propagation directions and the performance of the AIFM method, offering guidance for optimizing the model to achieve both high accuracy and efficiency.

We will examine the effect of varying the number of Ricker wave propagation directions $M$ on the AIFM model's particle detection accuracy and stability. Specifically, we compare the results for $M=10,20$, evaluating the influence of the number of propagation directions on the precision and robustness of the particle detection model.

\begin{figure}
	\centering
	\begin{minipage}{0.32\textwidth}
		\subfigure[Ground truth]{\includegraphics[width=0.95\columnwidth]{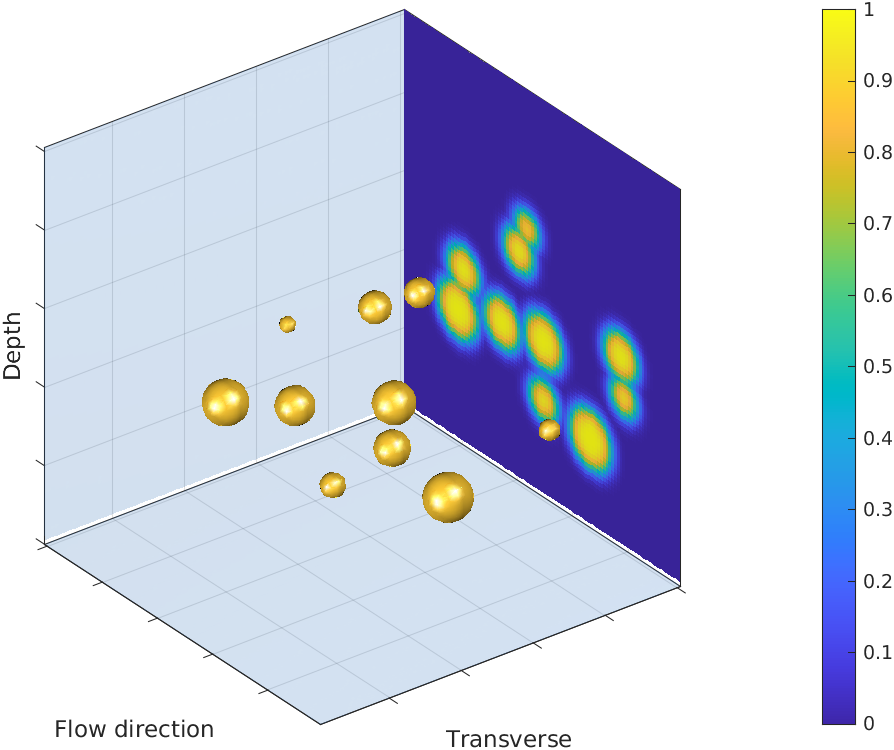}
		\label{fig_effect_num_directions_mask}}
	\end{minipage}
	\begin{minipage}{0.32\textwidth}
		\subfigure[10 directions]{\includegraphics[width=0.95\columnwidth]{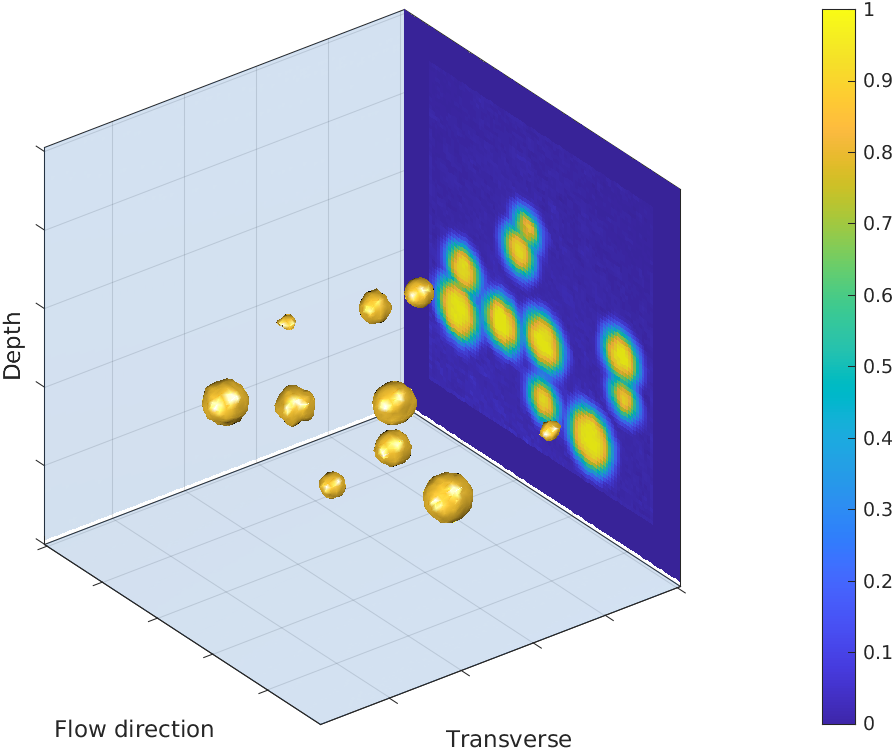}
		\label{fig_effect_num_directions_10}}
	\end{minipage}
	\begin{minipage}{0.32\textwidth}
		\subfigure[20 directions]{\includegraphics[width=0.95\columnwidth]{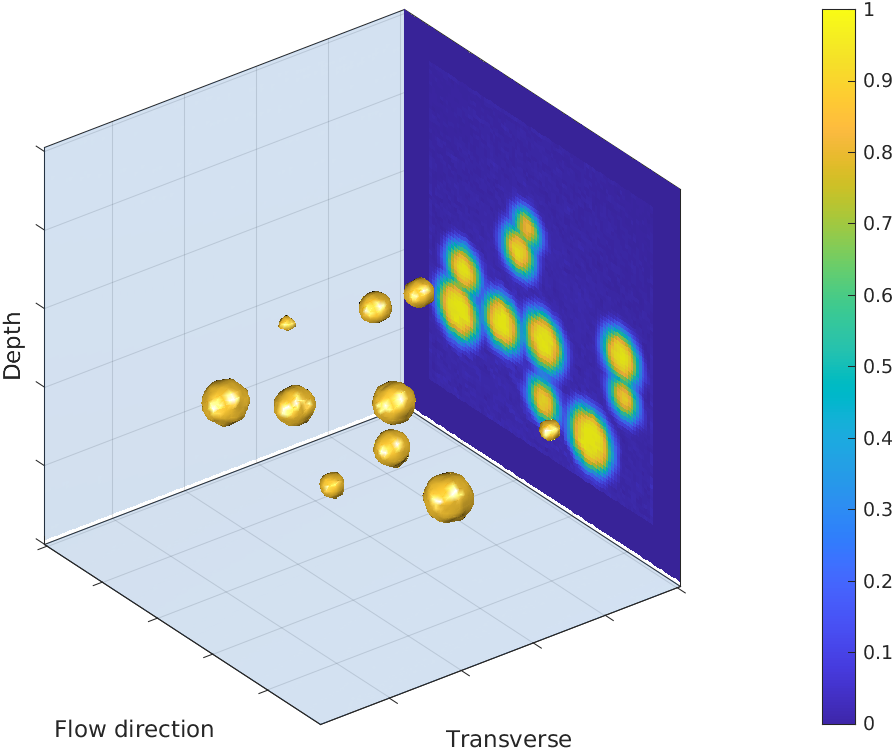}
		\label{fig_effect_num_directions_20}}
	\end{minipage}
	\caption{Effect of the number of directions of emitted waves. The first figure is the ground truth and the other two figures are the reconstructions of particles using the AIFM method with 10 and 20 directions of acoustic waves, respectively.}
	\label{fig_effect_num_directions}
\end{figure}

Next, we present the experimental results to evaluate the impact of the number of emission directions on the inversion accuracy of the AIFM model, as shown in \cref{fig_effect_num_directions}. \cref{fig_effect_num_directions_mask} illustrates the true particle mass distribution, which serves as a reference baseline for comparison. The reconstructed particle mass distributions using the AIFM model with 10 and 20 emission directions are shown in \cref{fig_effect_num_directions_10} and \cref{fig_effect_num_directions_20}, respectively. These results demonstrate how the number of emission directions affects the model's ability to accurately reconstruct the particle mass distribution.

The experimental results demonstrate that the AIFM model effectively reconstructs the particle mass distribution using both 10 and 20 emission directions, confirming its stability and adaptability. While the 20-direction configuration shows a slight improvement in accuracy compared to the 10-direction case, the marginal gain does not justify the significant increase in computational cost. Therefore, 10 emission directions are recommended as the optimal choice, striking a balance between computational efficiency and inversion accuracy. This configuration will be adopted in subsequent numerical experiments.

\subsubsection{Effect of the number of receivers}\label{subsubsubsec_receiverNumber}
The number of receivers is a critical factor influencing the accuracy of particle detection in the AIFM method, as increasing their number enhances inversion accuracy by providing more comprehensive data for reconstructing the particle mass distribution. A key advantage of the AIFM method is that its computational complexity remains unaffected by the number of receivers due to its parallel processing approach, which handles data from all receivers simultaneously. This allows for improved accuracy without significantly increasing computational costs. However, practical constraints such as measurement conditions, installation costs, and equipment availability often limit the number of receivers. In this experiment, we evaluate the impact of receiver quantity on detection accuracy using a setup where receivers are arranged on the three boundary surfaces of the channel and the water surface, reflecting typical real-world sensor placement configurations. This approach aims to balance computational accuracy with practical limitations.

We will evaluate three receiver configurations: $4 \times 101^2$, $4 \times 51^2$, and $4 \times 21^2$ receivers. These configurations reflect practical measurement conditions and consider factors like equipment costs and data feasibility. This comparison helps assess how receiver quantity affects particle detection precision in the AIFM model.

\begin{figure}
	\centering
	\begin{minipage}{0.24\textwidth}
		\subfigure[Ground truth]{\includegraphics[width=0.95\columnwidth]{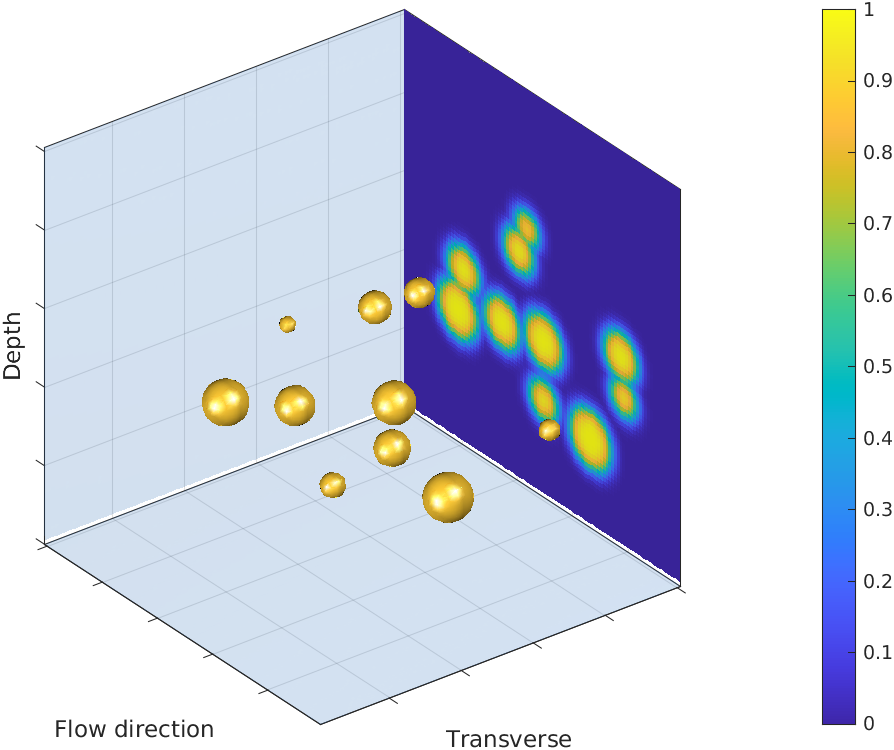}}
		\label{fig_effect_num_receivers_mask}
	\end{minipage}
	\begin{minipage}{0.24\textwidth}
		\subfigure[$N=4\times 101^2$]{\includegraphics[width=0.95\columnwidth]{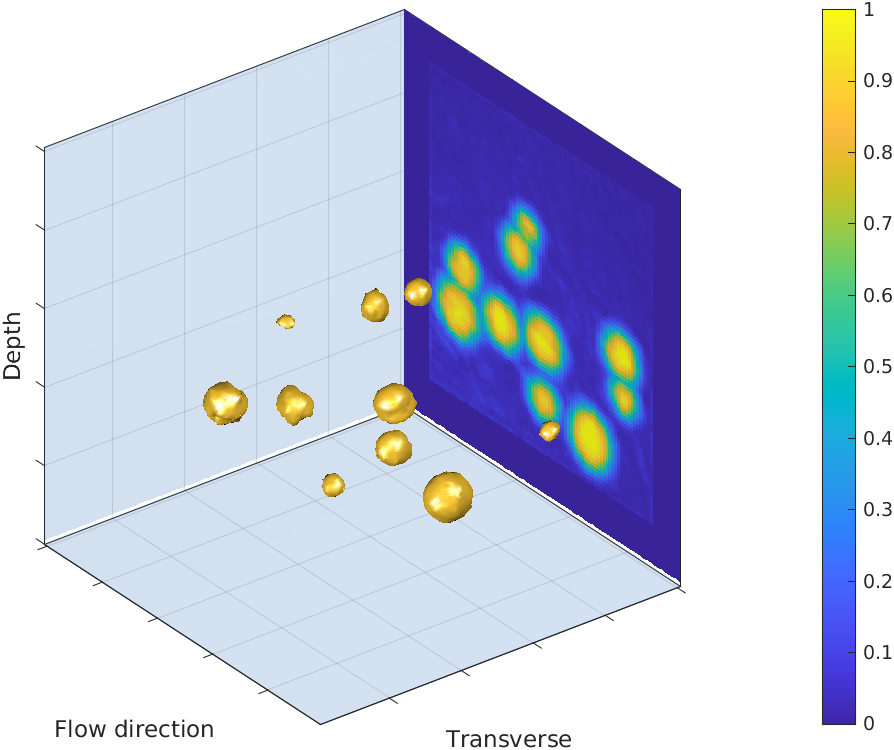}}
		\label{fig_effect_num_receivers_101}
	\end{minipage}
	\begin{minipage}{0.24\textwidth}
		\subfigure[$N=4\times 51^2$]{\includegraphics[width=0.95\columnwidth]{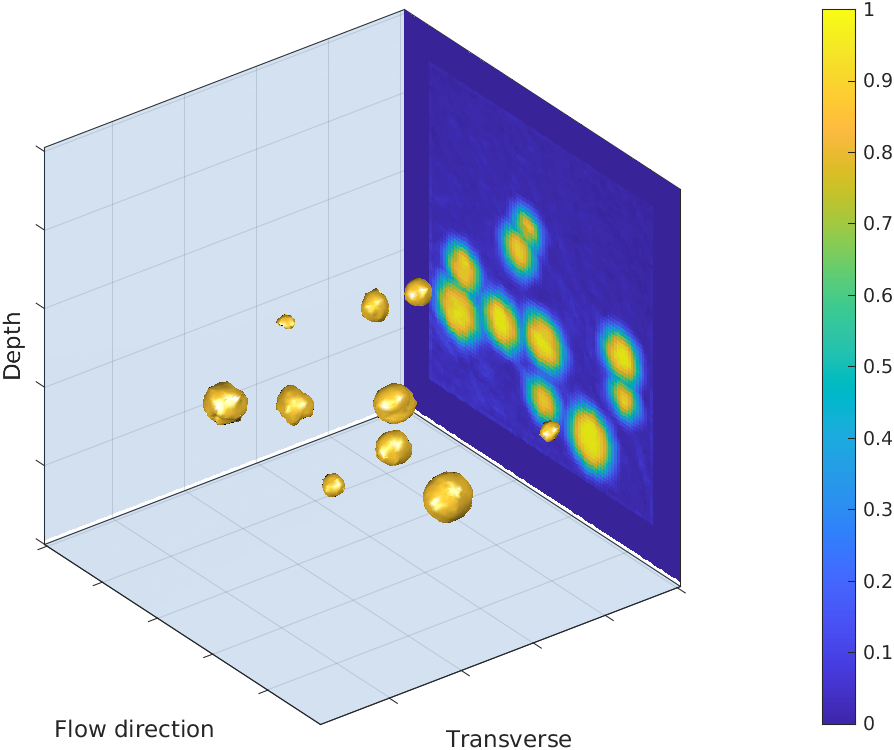}}
		\label{fig_effect_num_receivers_51}
	\end{minipage}
	\begin{minipage}{0.24\textwidth}
		\subfigure[$N=4\times 21^2$]{\includegraphics[width=0.95\columnwidth]{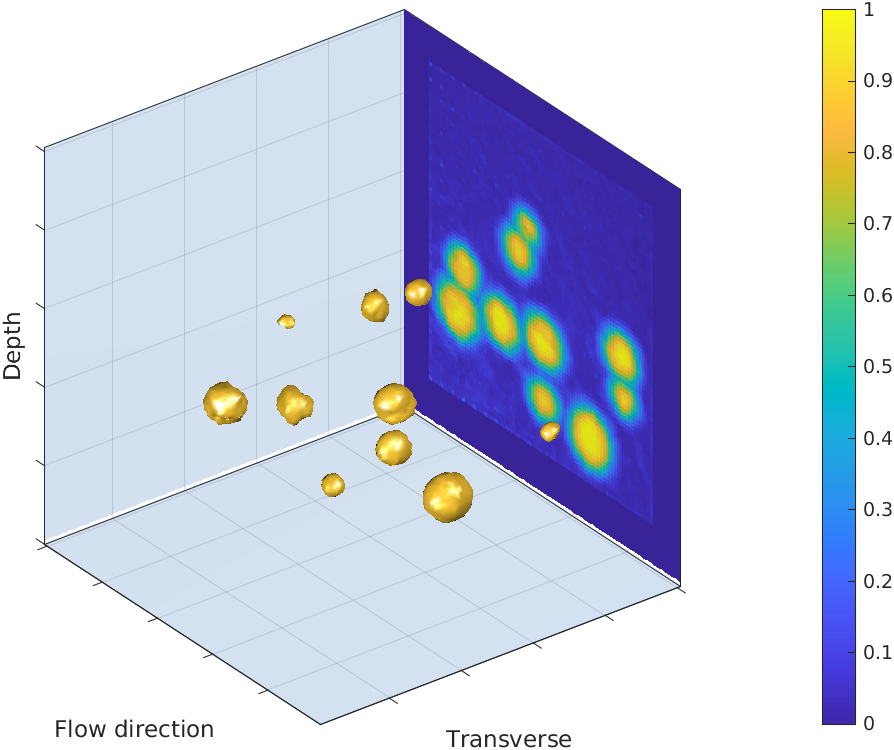}}
		\label{fig_effect_num_receivers_21}
	\end{minipage}
	\caption{Effect of the number of receivers. The first figure is the ground truth and the other three figures are the recoveries of particles using the AIFM method with $4\times 101^2, 4\times 51^2$ and $4\times 21^2$ receivers, respectively.}
	\label{fig_effect_num_receivers}
\end{figure}

The results, shown in \cref{fig_effect_num_receivers}, demonstrate that increasing the number of receivers improves particle detection accuracy, as more data enhances reconstruction. The true solution is provided in \cref{fig_effect_num_receivers_mask}, while the inversion results for $4 \times 101^2$, $4 \times 51^2$, and $4 \times 21^2$ receivers are displayed in \cref{fig_effect_num_receivers_101}, \cref{fig_effect_num_receivers_51}, and \cref{fig_effect_num_receivers_21}, respectively.

The results suggest that while more receivers provide additional data, the marginal improvement in accuracy diminishes after a certain point. In practical applications, using a reasonable number of receivers, such as $101^2$ per face, achieves the desired reconstruction accuracy without unnecessarily increasing computational load. This approach optimizes the trade-off between computational efficiency and detection accuracy, ensuring a balanced and effective solution.

\subsubsection{Effect of the layout of receivers}\label{subsubsubsec_receiverLayout}
\begin{figure}
	\centering
	\begin{minipage}{0.24\textwidth}
		\subfigure[Ground truth]{\includegraphics[width=0.95\columnwidth]{figures/effect_num_receivers/mask.png}}
		\label{fig_effect_num_receivers_layout_mask}
	\end{minipage}
	\begin{minipage}{0.24\textwidth}
		\subfigure[All-around layout]{\includegraphics[width=0.95\columnwidth]{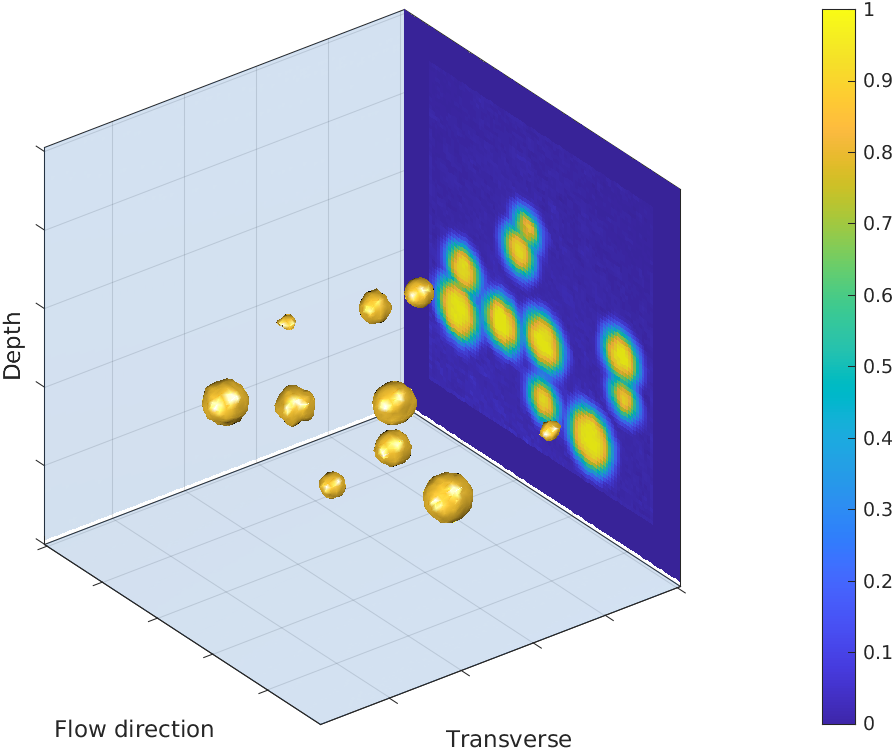}}
		\label{fig_effect_num_receivers_layout_101}
	\end{minipage}
	\begin{minipage}{0.24\textwidth}
		\subfigure[Canal walls and water surface layout]{\includegraphics[width=0.95\columnwidth]{figures/effect_num_receivers/inv_Ricker20k_recXoyYoz_dirs10_river_mask_inv.png}}
		\label{fig_effect_num_receivers_layout_51}
	\end{minipage}
	\begin{minipage}{0.24\textwidth}
		\subfigure[Canal sidewalls layout]{\includegraphics[width=0.95\columnwidth]{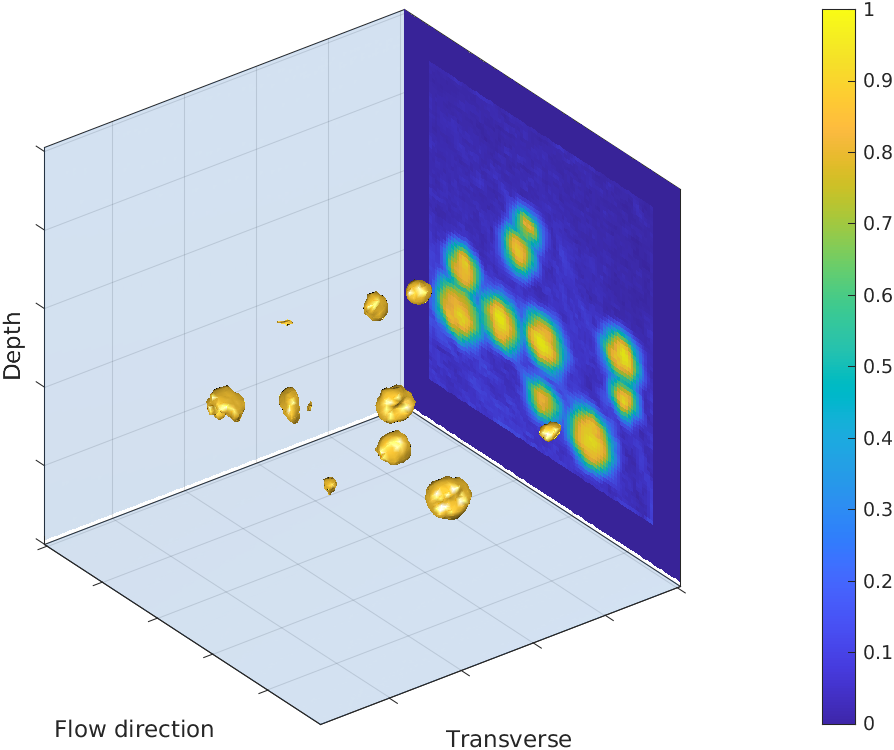}}
		\label{fig_effect_num_receivers_layout_21}
	\end{minipage}
	\caption{Effect of the layout of receivers. The first figure is the ground truth and the other three figures are the recoveries of particles using the AIFM method with different layouts of receivers.}
	\label{fig_effect_num_receivers_layout}
\end{figure}

To further demonstrate that the performance of the AIFM method is primarily influenced by the layout of the receivers rather than their number, we conduct an additional experiment. In this experiment, we fix the number of receivers on a single face of the domain at $101^2$ and vary the receiver layout. Specifically, we explore three different layouts: an all-around configuration, a configuration covering the canal walls and water surface, and a configuration along the canal sidewalls. The results are presented in \cref{fig_effect_num_receivers_layout}.

The receiver layout significantly impacts particle detection accuracy, with the highest accuracy achieved when receivers are distributed across all faces. This configuration provides the best approximation of the particle distribution. In contrast, layouts limited to the canal walls and water surface or only the canal sidewalls show reduced accuracy, particularly in detecting smaller particles, highlighting the importance of optimized receiver placement.

The comparison of results in \cref{fig_effect_num_receivers} and \cref{fig_effect_num_receivers_layout} reveals that the receiver layout has a greater impact on the performance of the AIFM method than the number of receivers. Therefore, optimizing the strategic arrangement of receivers should be prioritized to ensure the highest accuracy in particle detection. Once the layout is optimized, the number of receivers can be adjusted to balance accuracy and computational cost, ensuring the model remains both effective and efficient in practical applications.

\subsubsection{Effect of the number of particles}\label{subsubsubsec_particleNumber}
Most flow measurement techniques are limited by the number of particles that can be effectively detected. Each instrument has an upper bound on the number of particles it can handle accurately. Exceeding this limit typically results in reduced measurement accuracy and stability. In this section, we conduct a series of numerical experiments to explore the impact of particle quantity on the performance of the AIFM model's particle detection and to determine the optimal range of particle quantities for the AIFM model. The experimental results will help us understand how particle quantity influences the accuracy and stability of flow rate measurements and provide guidance for selecting the appropriate particle quantity in practical measurements.

\begin{figure}
	\centering
	\begin{minipage}{0.32\textwidth}
		\subfigure[Ground truth: 10 points]{\includegraphics[width=0.95\columnwidth]{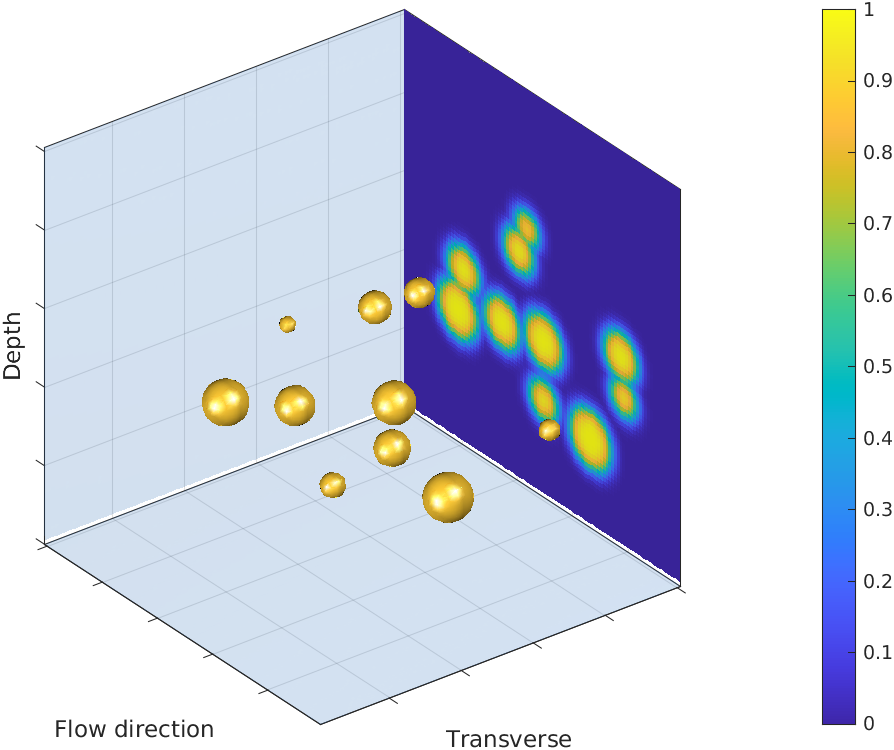}}
		\label{fig_effect_num_particles_true_10}
	\end{minipage}
	\begin{minipage}{0.32\textwidth}
		\subfigure[Ground truth: 50 points]{\includegraphics[width=0.95\columnwidth]{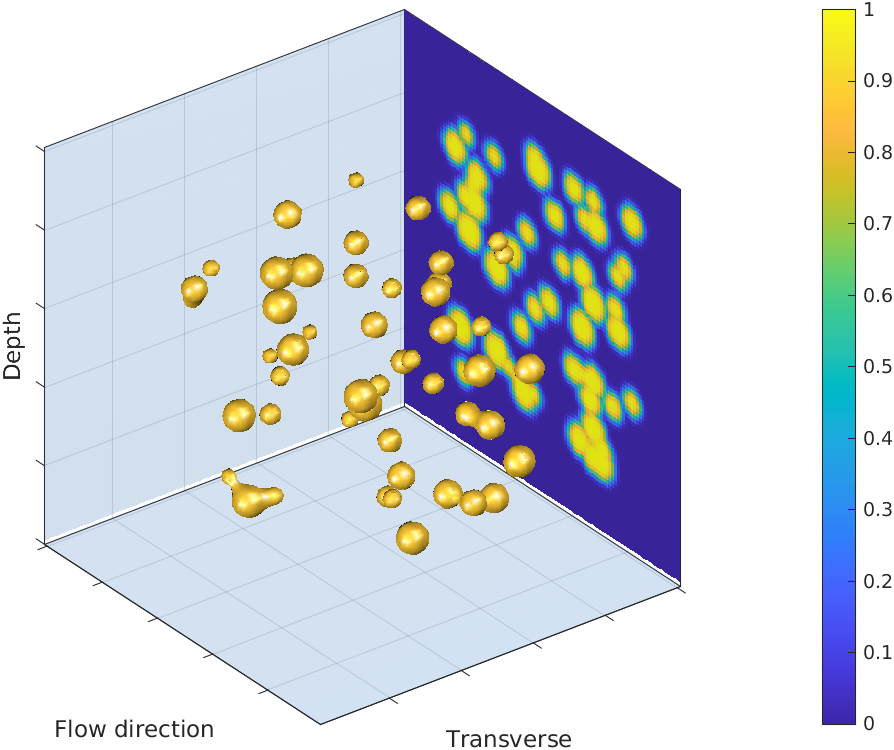}}
		\label{fig_effect_num_particles_true_50}
	\end{minipage}
	\begin{minipage}{0.32\textwidth}
		\subfigure[Ground truth: 200 points]{\includegraphics[width=0.95\columnwidth]{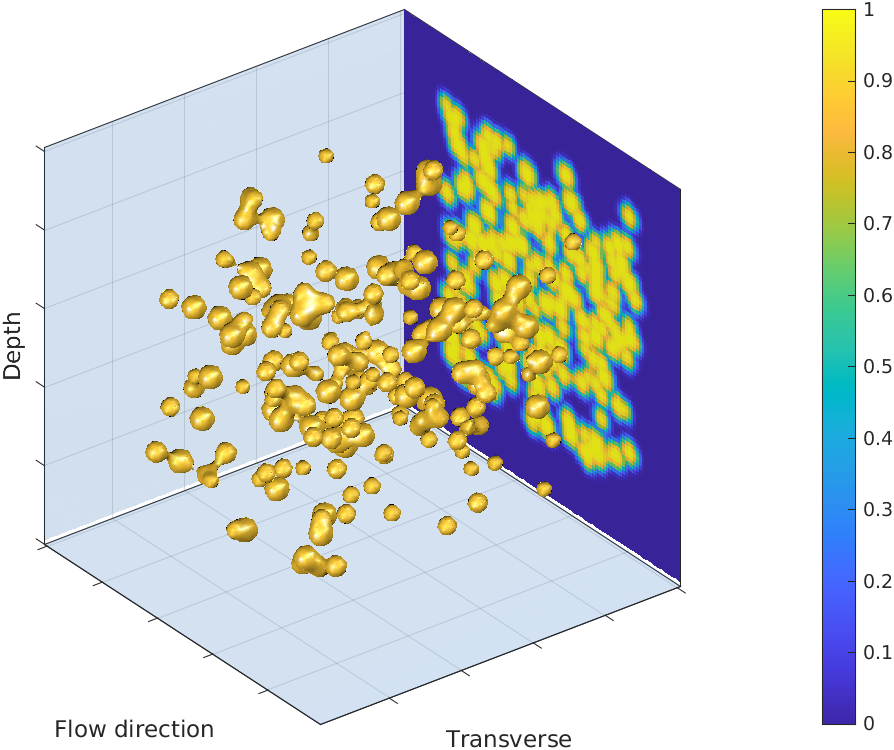}}
		\label{fig_effect_num_particles_true_200}
	\end{minipage}

	\begin{minipage}{0.32\textwidth}
		\subfigure[Recovery: 10 points]{\includegraphics[width=0.95\columnwidth]{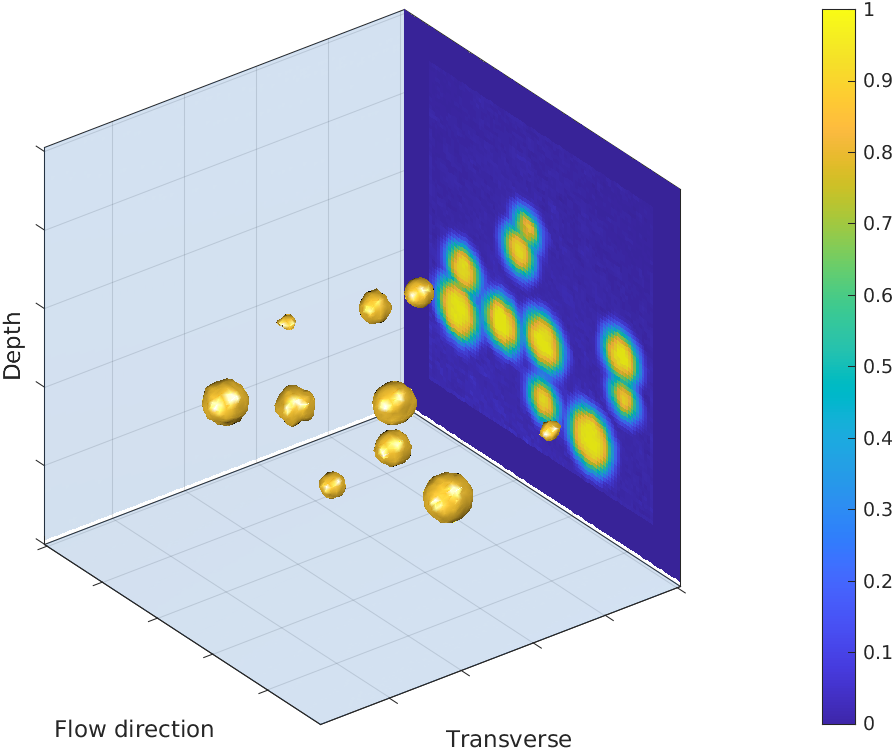}}
		\label{fig_effect_num_particles_inv_10}
	\end{minipage}
	\begin{minipage}{0.32\textwidth}
		\subfigure[Recovery: 50 points]{\includegraphics[width=0.95\columnwidth]{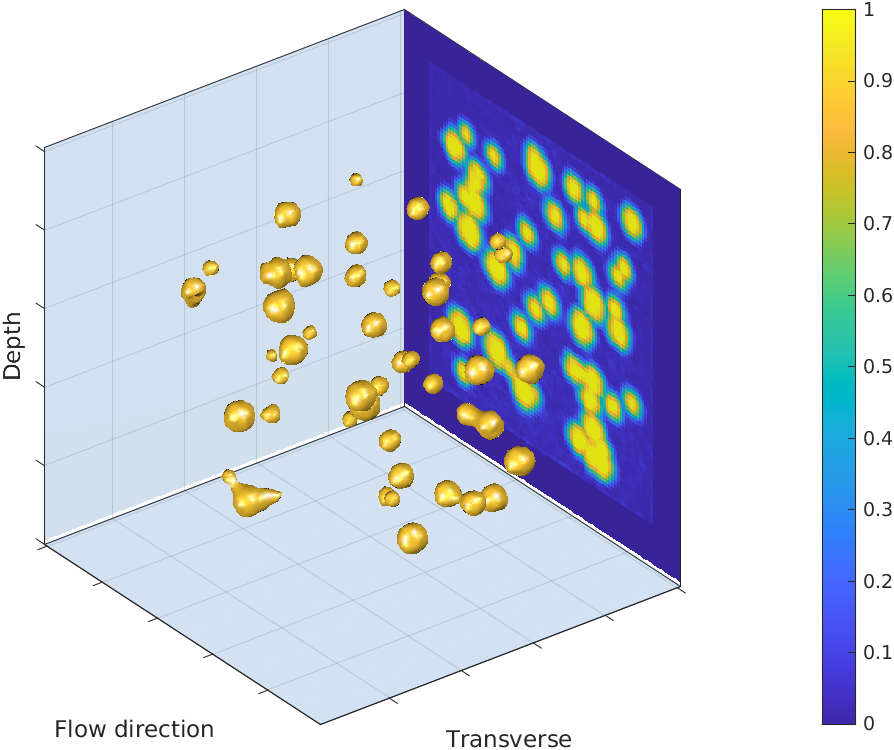}}
		\label{fig_effect_num_particles_inv_50}
	\end{minipage}
	\begin{minipage}{0.32\textwidth}
		\subfigure[Recovery: 200 points]{\includegraphics[width=0.95\columnwidth]{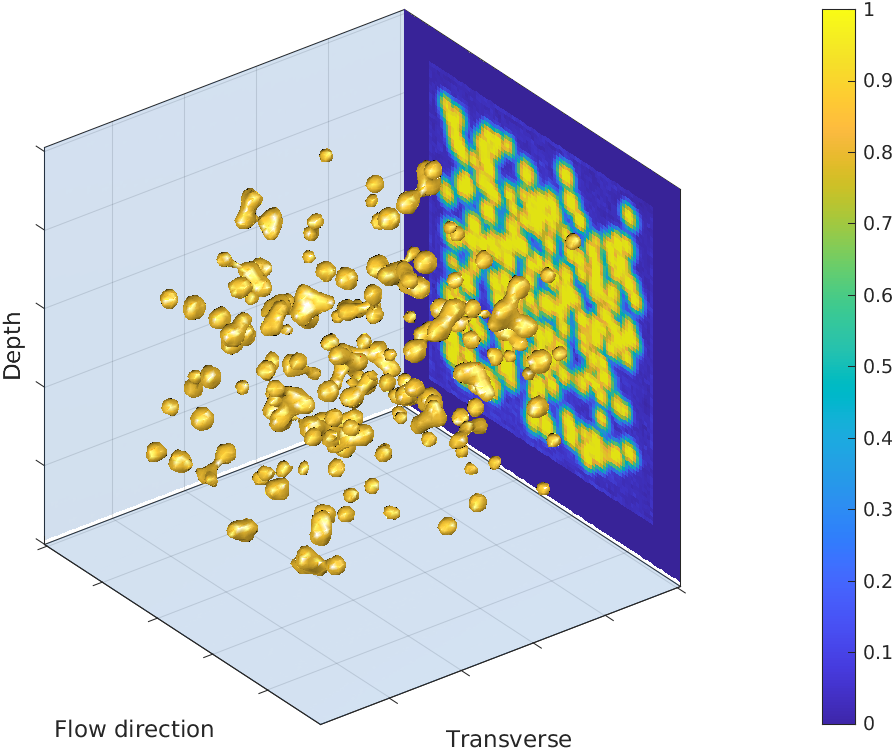}}
		\label{fig_effect_num_particles_inv_200}
	\end{minipage}
	\caption{Effect of the number of particles. The top row is the ground truths of particles with 10, 50, and 200 particles, respectively. The bottom row is the corresponding recoveries of particles using the AIFM method.}
	\label{fig_effect_num_particles}
\end{figure}

The experiment is designed with the receiver layout fixed to the all-around configuration, and the number of receivers set to $N=6\times 101^2$. We compare the particle detection performance for three different particle quantities: 10, 50, and 200 particles. The goal is to evaluate the performance of AIFM model at varying particle counts and determine the optimal range of particles for which the model maintains high accuracy in reconstructing the particle distribution. The results are presented in \cref{fig_effect_num_particles}. The top and bottom rows of \cref{fig_effect_num_particles} show the ground truth and the particle reconstructions, respectively.

The results demonstrate that the AIFM model achieves accurate particle detection and reconstruction for all tested particle quantities—10, 50, and 200 particles. This confirms the model’s reliability across a wide range of particle counts. Notably, even with 200 particles, the reconstruction maintains high precision, highlighting the model’s robustness for handling high particle densities. Further analysis reveals that the AIFM model maintains remarkable stability and accuracy even as particle counts vary, making it highly suitable for real-world applications where particle quantities fluctuate dynamically. This consistency ensures the model’s effectiveness in dynamic fluid systems, underscoring its broad applicability and practical value.

\subsection{Flow velocity field measurement}\label{subsec_flowMeasurement}
In this section, we would evaluate the accuracy and stability of the AIFM model in calculating flow velocity fields under varying flow conditions. Specifically, the scenarios include a constant flow field, a vortex flow field, and a T-junction flow field, each representing distinct practical flow environments. These tests will provide valuable insights for further optimizing the AIFM model, improving the measurement accuracy and stability of fluid velocity measurement equipment in complex environments.

The selection of these three flow field scenarios is designed to comprehensively assess the applicability and robustness of the AIFM model under different flow conditions. Each flow field corresponds to a different practical application scenario, which has significant engineering values. The constant flow field represents an idealized flow environment, commonly used for calibrating and validating fluid velocity measurement instruments. This scenario enables us to verify whether the AIFM model can accurately reconstruct constant flow velocities in a controlled environment. The vortex flow field, on the other hand, simulates more complex flow dynamics. This scenario is relevant for testing the performance of the AIFM model in more intricate flow environments, especially its ability to respond to and accurately capture flow variations. Finally, the T-junction flow field is designed to simulate real-world complex flow channel structures, such as pipes or riverbeds. The fluid velocity in such scenarios is influenced by boundary conditions and channel shapes, making it a challenging environment for the model. This experiment will evaluate the ability of the AIFM model to accurately calculate the velocity field in regions with significant velocity changes.

These experiments not only aim to assess the accuracy of the AIFM model across different flow fields but also to provide detailed performance data for its practical applications in diverse flow environments. The results will offer a deeper understanding of the model's adaptability and robustness, further enhancing the utility of fluid velocity measurement devices in engineering applications.

\subsubsection{Experimental setup}
The numerical experiments build upon the experimental setup described in \cref{subsec_particleDetection}, maintaining consistent conditions by using the same number of sound wave sources and receivers. This ensures a fair comparison of the performance of different flow velocity field computation models under uniform experimental settings.

To optimize the performance of the Farneb\"{a}ck optical flow method, we specifically select a particle density of $200$ particles per unit volume. This choice is based not only on the optimal performance of the Farneb\"{a}ck optical flow algorithm for tracking particle movement, but also on balancing the computational errors between the particle detection model and the flow velocity field computation model. If the number of particles is too low, the accuracy of the flow field calculation may be compromised; conversely, a higher particle count, although improving the flow field calculation precision, could negatively affect the accuracy of the particle detection model. Therefore, the choice of $200$ particles per unit volume represents a balanced compromise, ensuring optimal performance for both the particle detection and flow velocity field computation models.

To capture dynamic variations in the flow velocity field while ensuring computational precision, we collect sound wave data at $10$ distinct time points, adjusting the time intervals based on the overall flow speed. For flow velocities below \SI{10}{\meter\per\second}, a longer interval of \SI{1}{\second} is used, whereas for speeds exceeding \SI{10}{\meter\per\second}, a shorter interval of \SI{0.5}{\second} is selected. This approach prevents rapid particle motion from causing inaccuracies in the flow field calculations.

Selecting an appropriate time interval is crucial for maintaining the AIFM model's accuracy. Too long an interval may fail to capture rapid particle movements in high-speed regions, increasing errors, while too short an interval can introduce noise or instability. By dynamically adjusting the interval based on flow velocity, we optimize measurement precision and computational stability.

\subsubsection{Evaluation metrics}
To comprehensively assess the accuracy and stability of the AIFM model’s flow field computation, we will employ the following four types of relative errors as quantitative metrics, to evaluate the performance under various measurement conditions:
\begin{enumerate} 
	\item RE1: The relative $L^2$ error of the reconstructed velocity field across the entire domain. 
	\item RE2: The relative $L^2$ error of the reconstructed velocity field restricted to particle regions. 
	\item RE3: The relative error of the average velocity vector projected along the flow direction over the entire domain.
	\item RE4: The relative error of the velocity vector projected along the flow direction within particle regions. 
\end{enumerate}

Due to the higher accuracy of the Farneb\"{a}ck optical flow method near particles, we chose relative errors RE2 and RE4 as the primary indicators for evaluating the accuracy and stability of the AIFM model’s flow field computation. This localized accuracy characteristic of the Farneb\"{a}ck method, while providing excellent flow field reconstruction capabilities in particle-dense regions, also exposes a common issue: in areas distant from the particles, the accuracy of the flow field computation is lower. This phenomenon is a shared limitation across all acoustic-based flow measurement methods. Specifically, in regions with fewer particles, the signal strength of the interaction between the sound waves and the fluid diminishes, leading to increased inaccuracies in the measurements \cite{liuComparisonOpticalFlow2015,jiHighresolutionVelocityDetermination2024,glombOpticalFlowbasedMethod2017}. The root cause of this issue lies in the difficulty of effectively transmitting information in areas devoid of particles, thus impacting the overall accuracy of the flow field computation model.

Therefore, the selection of these four quantitative indicators takes into account not only the high-precision measurements in particle-dense regions but also the computational errors in regions farther from the particles. With these two indicators, we can comprehensively assess the performance of the AIFM model under different conditions, particularly its stability and reliability in overall flow field reconstruction. Although this issue is common in acoustic measurement technologies, we can still reduce such regional errors and enhance the model’s overall accuracy by improving algorithms, increasing sensor density, and optimizing data processing methods.

Next, we will provide a detailed analysis of the AIFM model’s flow field computation in distinct flow field scenarios as described above, evaluating its accuracy and stability under varying flow conditions.
\subsubsection{Constant flow}
In this subsection, we investigate the inversion of a three-dimensional constant flow velocity field. The constant flow field is a simple yet fundamental scenario widely used in the calibration and validation of fluid velocity measurement devices, pipeline velocity measurement, and the measurement of flow velocity in straight river channels. To ensure generality, the flow velocity is set to a constant vector $(0,\SI{1}{\meter\per\second},0)$, which is independent of both temporal and spatial coordinates. In this flow field, we evaluate the accuracy and stability of the AIFM model’s flow velocity computation model to verify its applicability in an ideal flow scenario.

\begin{figure}
	\centering
	\begin{minipage}{0.4\textwidth}
		\subfigure[Ground truth]{\includegraphics[width=0.95\columnwidth]{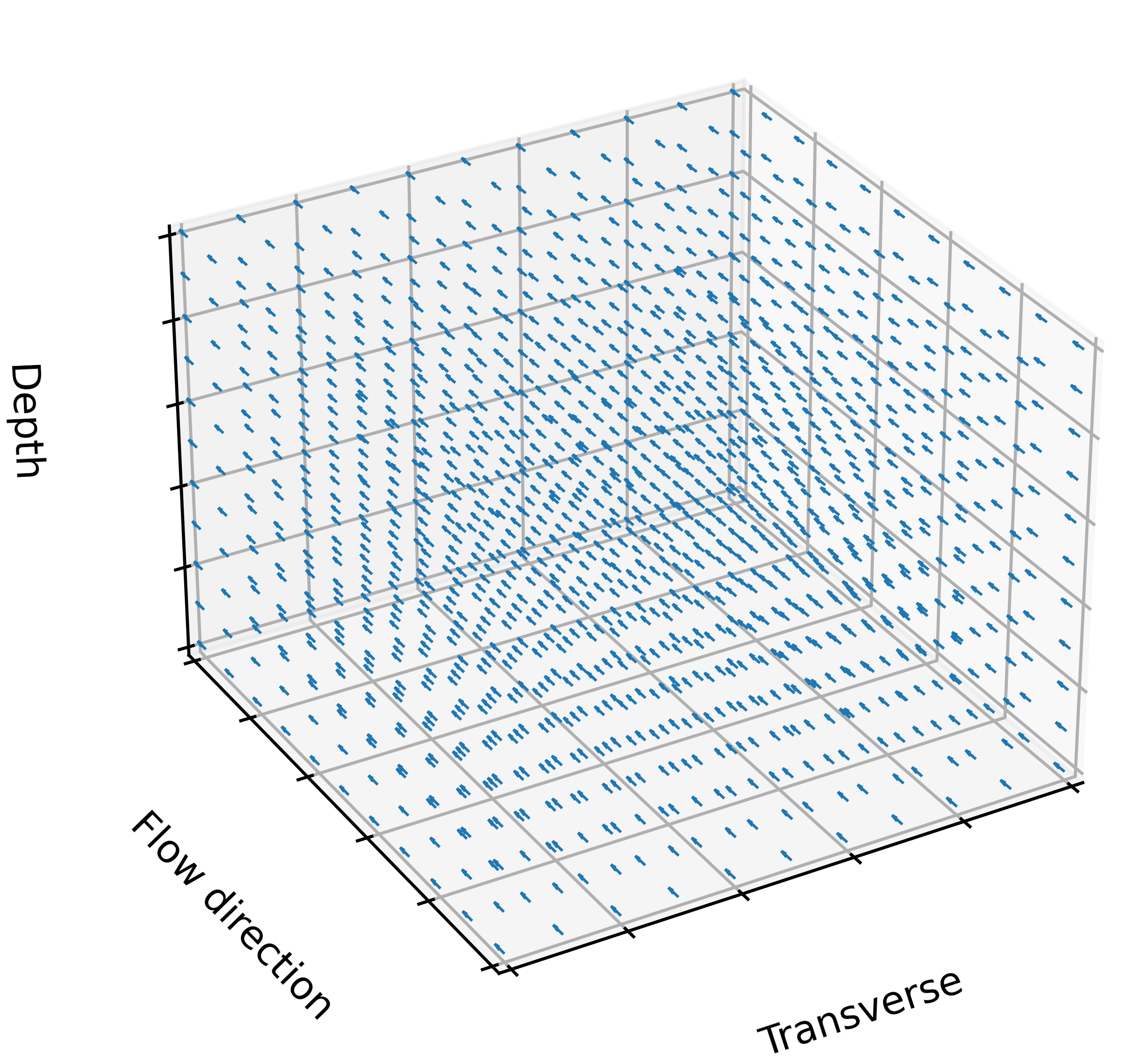}
		\label{fig_constantFlow_groundTruth}}
	\end{minipage}
	\begin{minipage}{0.4\textwidth}
		\subfigure[Reconstructed flow field]{\includegraphics[width=0.95\columnwidth]{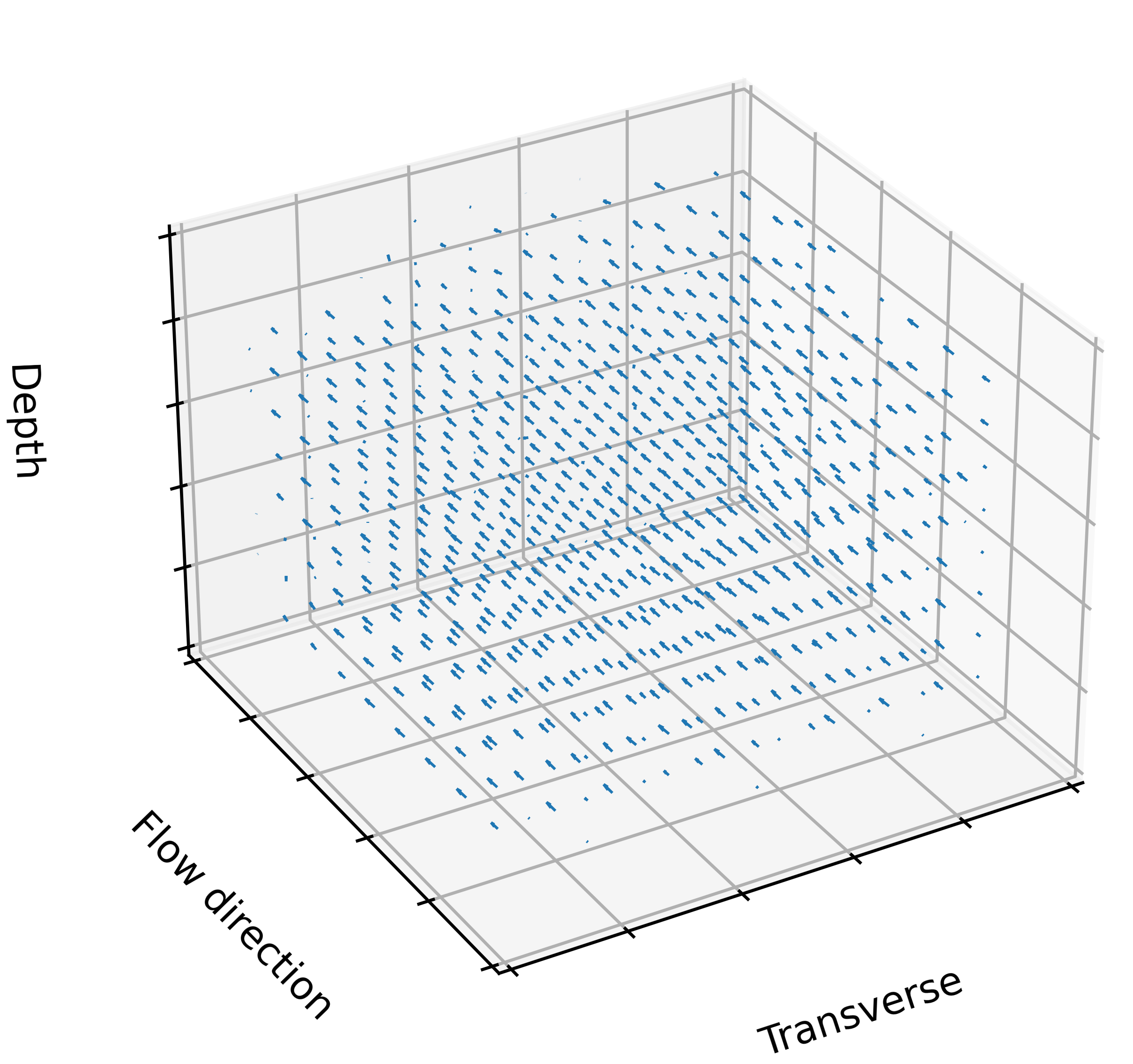}
		\label{fig_constantFlow_measured}}
	\end{minipage}
	\caption{Scenario of the constant flow. The left figure is the ground truth and the right figure is the measured flow field using the optical flow method.}
	\label{fig_constantFlow}
\end{figure}

The AIFM model accurately reconstructs the characteristics of the real flow field, as demonstrated by the inversion results in \cref{fig_constantFlow}. \cref{fig_constantFlow_groundTruth} shows the real flow field, while \cref{fig_constantFlow_measured} presents the inverted flow field. The relative error between the two, detailed in \cref{table_constantFlow}, remains within \SI{10}{\percent} across the entire measurement area and below \SI{3}{\percent} for the average flow velocity vector.

These results highlight the AIFM model’s high accuracy and stability in steady flow fields, making it suitable for applications such as calibrating and validating fluid velocity measurement devices. The model’s ability to precisely reconstruct flow fields underlines its practical value in scenarios requiring reliable and accurate flow analysis.

\begin{table}
	\centering
	\caption{Relative errors in the constant flow scenario.}
	\label{table_constantFlow}
	\begin{tabular}{ccccc}
		\toprule
		Relative error & RE1 & RE2 & RE3 & RE4 \\\midrule
		Value          & \SI{7.34}{\percent} & \SI{5.33}{\percent} & \SI{2.34}{\percent} & \SI{1.95}{\percent} \\
		\bottomrule
	\end{tabular}
\end{table}

\subsubsection{Taylor-Green vortex flow}
In this subsection of the numerical experiments, we consider the inversion of the Taylor-Green vortex flow velocity field extended into three-dimensional space. The Taylor-Green vortex is a relatively complex flow field scenario that is applicable to the use of fluid velocity measurement devices in more intricate flow environments, such as in turbulence velocity measurements, vortex flow velocity measurements, and other similar applications. Due to the complex structure of the flow field, the Taylor-Green vortex introduces significant challenges for fluid velocity measurement technologies, particularly in regions with pronounced vortex structures and large velocity gradients, where measurement accuracy is often difficult to ensure. In this numerical experiment, we will analyze how the AIFM model handles this complex flow field and evaluate its inversion accuracy and stability under these conditions. By testing in such complex flow fields, we aim to gain deeper insights into the advantages and limitations of the AIFM model in practical applications, especially in flow fields with vortex characteristics and to determine whether the model can maintain high accuracy and robustness under such challenging conditions.

Typically, there is no analytical solution for the Taylor-Green vortex in three dimensions. For simplicity, we assume the vortex field is independent of the $x_3$-axis, which is given by
\begin{equation}
	\left(3\sin(\pi x_1)\cos(\pi x_2),-3\cos(\pi x_1)\sin(\pi x_2),0\right).
\end{equation}
Through this simplified setup, we are able to obtain a flow field with rotational symmetry that includes vortex characteristics. Although the model is relatively simplified, it still serves as a good representation of common vortex flow structures found in real-world fluid dynamics and is conducive to analysis and experimental validation.

\begin{figure}
	\centering
	\begin{minipage}{0.4\textwidth}
		\subfigure[Ground truth]{\includegraphics[width=0.95\columnwidth]{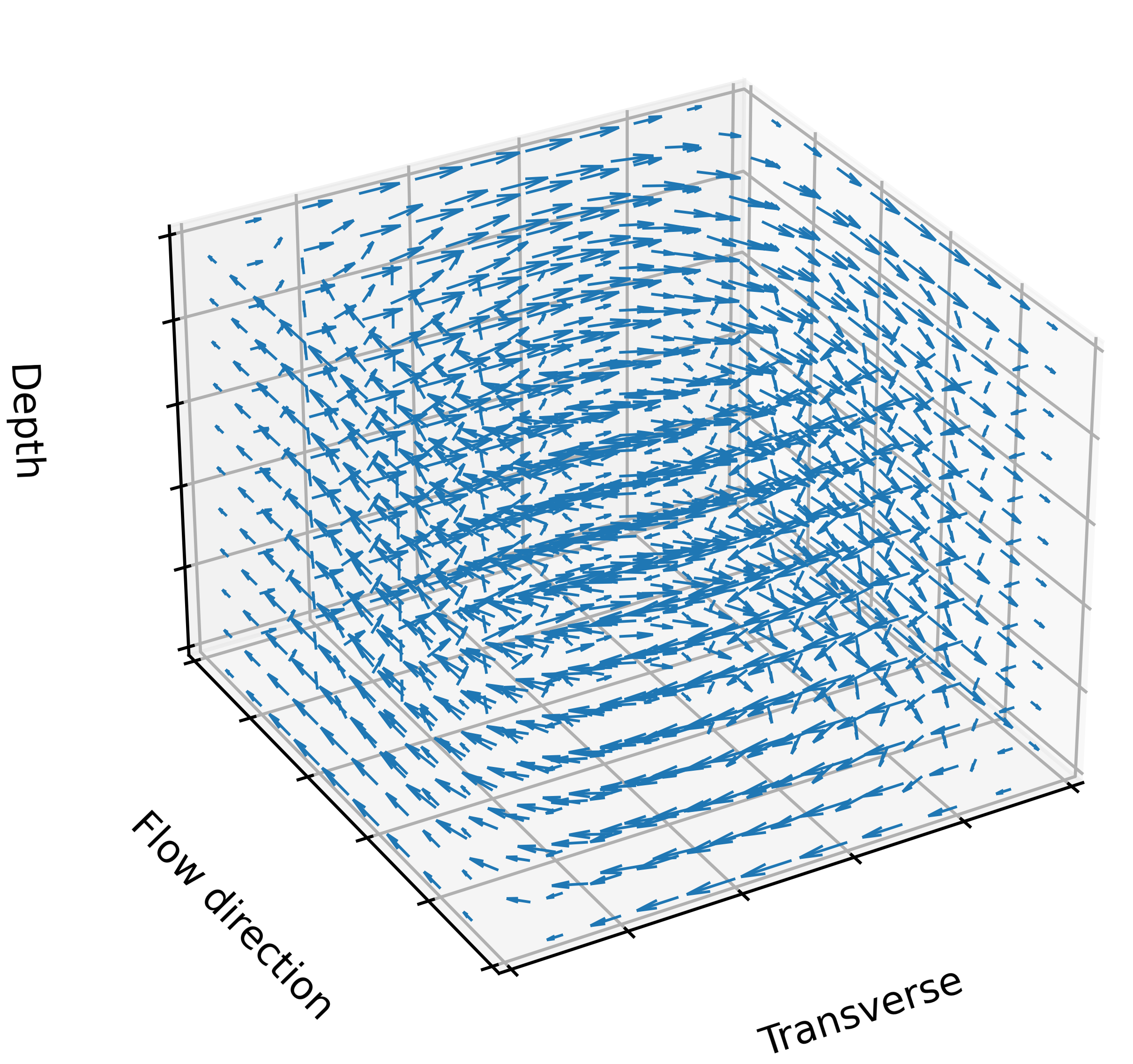}
		\label{fig_TaylorGreenVortex_groundTruth}}
	\end{minipage}
	\begin{minipage}{0.4\textwidth}
		\subfigure[Reconstructed flow field]{\includegraphics[width=0.95\columnwidth]{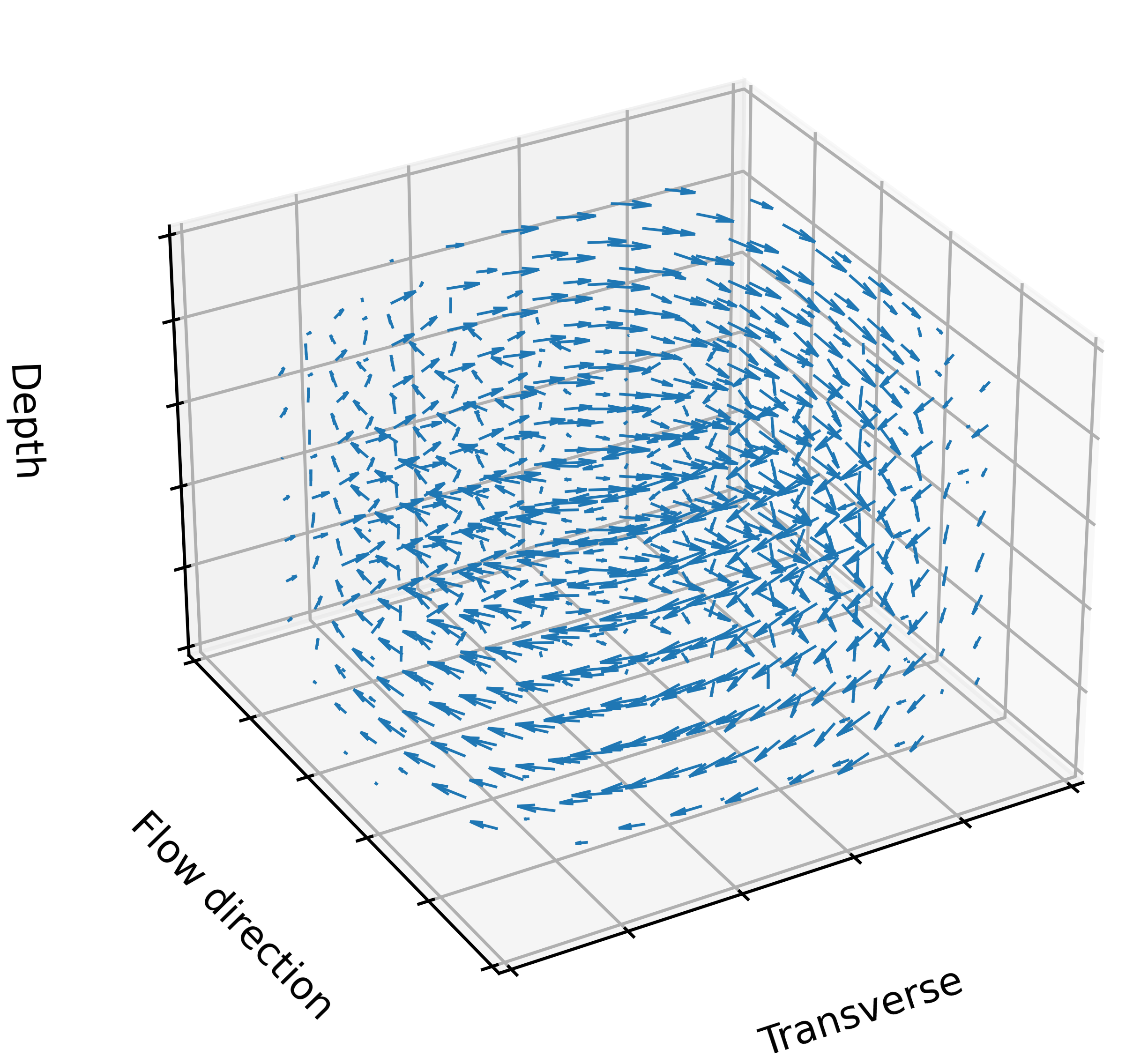}
		\label{fig_TaylorGreenVortex_measured}}
	\end{minipage}
	\caption{Scenario of Taylor-Green vortex. The left figure is the ground truth and the right figure is the measured flow field using the AIFM method.}
	\label{fig_TaylorGreenVortex}
\end{figure}

The AIFM model demonstrates high accuracy in reconstructing the vortex flow field, as shown in \cref{fig_TaylorGreenVortex}. \cref{fig_TaylorGreenVortex_groundTruth} displays the ground truth flow field, while \cref{fig_TaylorGreenVortex_measured} presents the inverted flow field. The relative errors (RE1 and RE2) between the inverted and ground truth flow fields, provided in \cref{table_TaylorGreenVortex}, remain below \SI{25}{\percent}, which is acceptable given the high complexity of the vortex flow.

These results highlight the precision and stability of the AIFM model in complex flow environments, effectively capturing the characteristics of the vortex flow field. This performance underscores the model’s suitability for applications in turbulent and vortex-dominated flows, where it can enhance measurement accuracy and reliability under challenging conditions.

\begin{table}
	\centering
	\caption{Relative errors for the Taylor-Green vortex.}
	\label{table_TaylorGreenVortex}
	\begin{tabular}{ccccc}
		\toprule
		Relative error & RE1 & RE2 & RE3 & RE4 \\\midrule
		Value          & \SI{24.25}{\percent} & \SI{20.00}{\percent} & \textendash& \textendash\\
		\bottomrule
	\end{tabular}
\end{table}

\subsubsection{T-type junction flow}
In this section, we explore the inversion of the flow field in a three-dimensional T-junction flow. The T-junction flow field is a complex flow scenario commonly encountered in many practical engineering applications, particularly at river confluences, pipeline junctions, and water pump inlets. These regions often present irregularities and instabilities in the flow field, posing significant challenges for flow measurement techniques. For example, at river confluences, the water flow may generate complex vortices and turbulence due to interactions with the banks and changes in flow patterns; in pipeline junctions, the flow can create localized vortex structures due to the branching of pipes. Given the challenges these environments present, evaluating the performance of the AIFM model in such flow fields is crucial. This section conducts numerical experiments to assess the accuracy, stability, and applicability of the AIFM model in T-junction flow fields, providing valuable insights for flow measurement in real-world engineering scenarios.

\begin{figure}
	\centering
	\includegraphics[width=0.5\columnwidth]{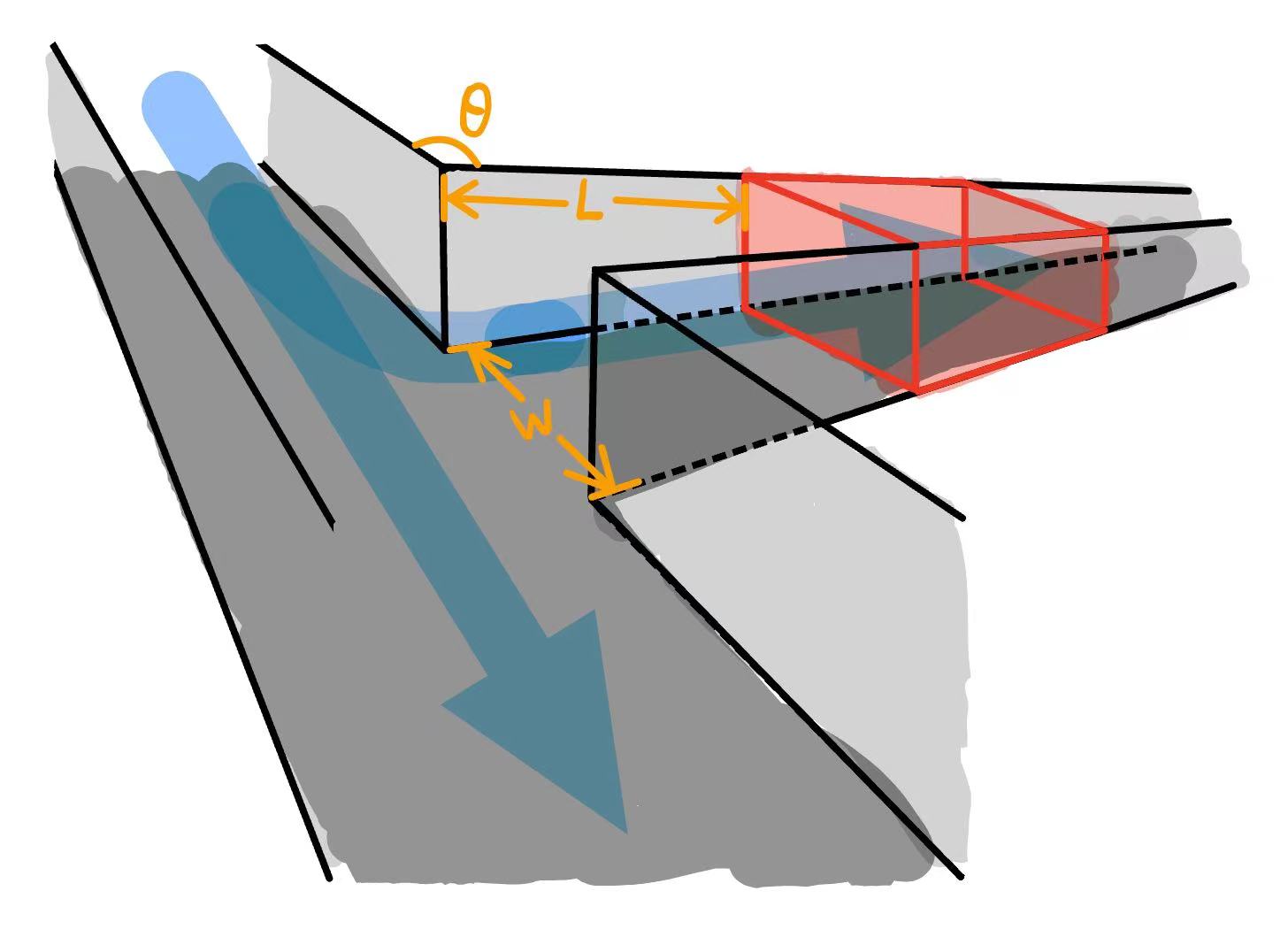}
	\caption{Setting of the T-type junction flow.}
	\label{fig_TtypeJunctionFlow_setting}
\end{figure}

To effectively simulate the T-type junction flow field, we use the commercial computational fluid dynamics software ANSYS Fluent 2024 R1, which provides robust numerical simulation capabilities for precise modeling of the flow field. ANSYS Fluent is widely used in the field of fluid dynamics, capable of handling complex flow field calculations, and is suitable for scenarios such as flow velocity measurement and turbulence analysis. In this experiment, we construct a typical T-type junction geometry model, as shown in \cref{fig_TtypeJunctionFlow_setting}. To simulate a realistic flow field, we set the inlet velocity to \SI{20}{\meter\per\second} and apply zero-pressure outlet conditions at the two exits, simulating common flow scenarios at pipe or river junctions. The main channel is \SI{100}{\meter} long and \SI{5}{\meter} wide, while the branch channel is \SI{50}{\meter} long and has a width of \SI{W}{\meter}. The depth of both channels is set to \SI{1}{\meter}, and the angle between the main and branch channels is $\theta \in (0, \pi)$.

The measurement region in this experiment is defined as \qtyproduct{W x 1 x 1}{\meter}, with the starting point located \SI{L}{\meter} from the T-type flow channel junction, as shown by the red cube in \cref{fig_TtypeJunctionFlow_setting}. To ensure inversion accuracy, we use the same experimental parameters as in \cref{subsec_particleDetection}, including 10 emitted acoustic wave fields and $6 \times 101^2$ receivers uniformly distributed across the boundaries, providing comprehensive coverage for high-precision flow field inversion. Additionally, we set the number of moving particles per unit volume to 200 to balance computational efficiency and particle detection accuracy, while collecting acoustic wave data at 10 evenly distributed time instances to effectively capture the dynamic evolution of the flow, particularly in regions with unstable velocity fields, ensuring reliable and accurate inversion results.

\begin{figure}
	\centering
	\begin{minipage}{0.3\textwidth}
		\subfigure[$\theta=\SI{60}{\degree}$]{\includegraphics[width=0.95\columnwidth]{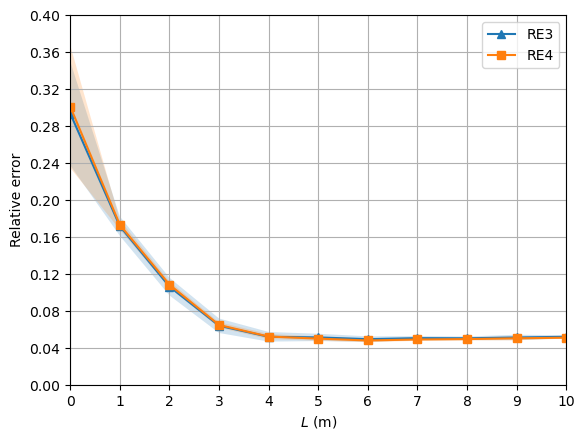}}
	\end{minipage}
	\begin{minipage}{0.3\textwidth}
		\subfigure[$\theta=\SI{90}{\degree}$]{\includegraphics[width=0.95\columnwidth]{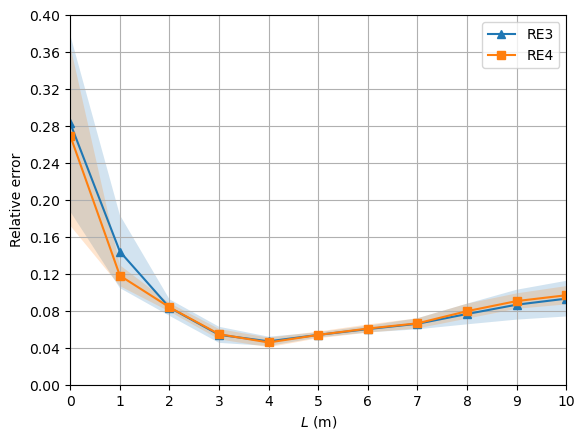}}
	\end{minipage}
	\begin{minipage}{0.3\textwidth}
		\subfigure[$\theta=\SI{120}{\degree}$]{\includegraphics[width=0.95\columnwidth]{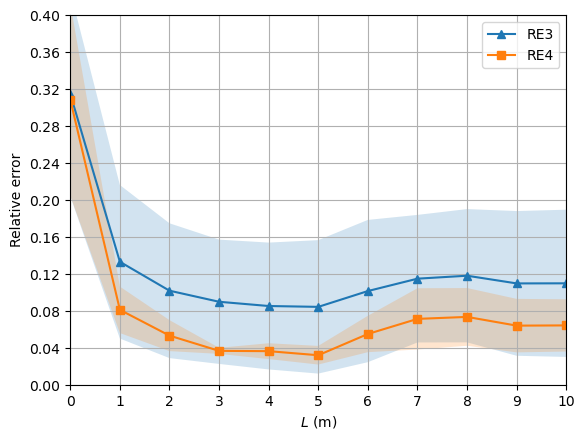}}
	\end{minipage}
	
	\caption{Relative errors of the AIFM method for the T-type junction flow with the fixed width $W=0.5$ and $\theta=\qtylist{60;90;120}{\degree}$. }
	\label{fig_TtypeJunctionFlow_errors_single}
\end{figure}

\begin{figure}
	\centering
	\includegraphics[width=0.3\columnwidth]{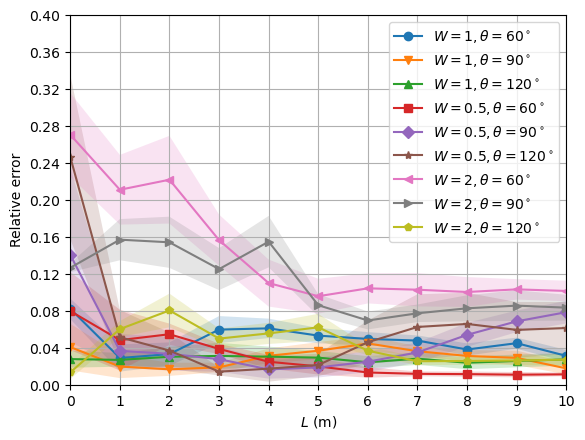}
	
	\caption{Relative errors of the AIFM model for several T-type junction flows with different $L$, $W$ and $\theta$.}
	\label{fig_TtypeJunctionFlow_errors_multiple}
\end{figure}

We will investigate how the AIFM model’s flow field inversion results vary with parameters $\theta,W$ and $L$, focusing on precision and stability in T-type junction flow fields. Accurate inversion is crucial for practical engineering applications, especially in complex environments where dynamic flow changes directly impact measurement reliability. We use relative errors RE3 and RE4 as primary evaluation metrics to comprehensively assess the model’s performance, aiming to quantify its adaptability and consistency across diverse flow scenarios. The first experiment fixes the branch channel width at $W=0.5$ and angle $\theta$ at \SIlist{60;90;120}{\degree}, varying the distance $L$ in the measurement region, with results shown in \cref{fig_TtypeJunctionFlow_errors_single}. The second experiment compares RE4 under different $W,\theta$ and $L$ values, as shown in \cref{fig_TtypeJunctionFlow_errors_multiple}. The shaded regions indicating the standard deviation of relative errors, reflecting measurement stability. These experiments highlight adaptability and reliability of the model in complex flow environments.

The results presented in \cref{fig_TtypeJunctionFlow_errors_single,fig_TtypeJunctionFlow_errors_multiple} reveal the following key observations.
\begin{enumerate}
	\item In the setup shown in \cref{fig_TtypeJunctionFlow_errors_single}, with the branch channel width fixed at $W=0.5$ and the angle $\theta$ set to \SIlist{60;90;120}{\degree}, the relative errors of the inversion results generally stay within \SI{15}{\percent}. Smaller values of $L$ result in slightly higher errors and decreased stability, likely due to the complex flow at the T-junction affecting the measurement region. However, as $L$ increases and the measurement area moves farther from the junction, the stability improves, with relative errors remaining below \SI{10}{\percent}. Moreover, the relative errors RE3 and RE4 show consistent trends, indicating that the AIFM model maintains high accuracy and stability both near the particles and in particle-free areas. This demonstrates the AIFM model's robustness and its ability to effectively handle flow field variations in complex environments, proving its wide applicability in real-world scenarios.
	\item In the experiment shown in \cref{fig_TtypeJunctionFlow_errors_multiple}, the relative error RE4 stays below \SI{30}{\percent} in all cases, with most values under \SI{20}{\percent}, indicating stable inversion of the flow field. Narrower branch channels lead to lower errors due to the more regular flow and less turbulence, making accurate reconstruction easier for the AIFM model. Changes in the angle $\theta$ cause minimal fluctuation in errors, highlighting the model's stability and adaptability to different geometries. However, in more complex situations—such as wider branches, smaller angles, or shorter distances from the junction—the relative error increases due to more irregular flow and stronger turbulence, though errors remain below \SI{30}{\percent}. This shows the ability of the AIFM model to capture overall flow characteristics, even in challenging environments. Overall, the AIFM model maintains high accuracy and stability across various conditions, demonstrating its strong applicability in real-world engineering scenarios.
\end{enumerate}

Overall, the AIFM model demonstrates high accuracy and stability in inversions of the flow velocity field in the T-junction. While errors increase in some complex flow conditions, they remain within acceptable limits, indicating that the model can reliably invert the flow field in most scenarios. The increased errors are primarily due to strong turbulence at the junction, which causes irregularities in the flow, especially in regions with high velocity and rapid particle movement. Despite these challenges, the AIFM model maintains a high level of stability, with errors still within an acceptable range. In conclusion, the AIFM model exhibits strong precision, stability, and robustness, making it well-suited for a wide range of real-world applications where it can adapt to various geometries and flow conditions while maintaining consistent accuracy. This highlights its potential and practicality for complex flow velocity measurements in engineering applications.

\section{Discussion}\label{sec_discussion}
This section presents a comprehensive evaluation of the AIFM model’s performance across different experimental configurations. The numerical results reveal key insights into how various parameters affect the accuracy, stability, and computational efficiency of the model. First, we investigated the influence of the number of propagation directions on the model’s particle detection accuracy and stability. Our findings indicate that while increasing the number of emission directions from 10 to 20 improves accuracy, the marginal gains do not justify the additional computational cost. Therefore, 10 emission directions provide an optimal balance between performance and efficiency. Similarly, when evaluating the effect of the number of receivers, we observed that adding more receivers enhances detection accuracy. However, the improvements become less significant beyond a certain number of receivers. A configuration of $101^2$ receivers per face strikes the best trade-off between computational load and detection precision. Furthermore, the receiver layout proved to be a more crucial factor than the receiver count, with the all-around distribution providing the highest accuracy. The configuration significantly impacts the model’s ability to detect and reconstruct particle mass distributions, particularly for smaller particles. Finally, the model was tested across varying particle quantities, ranging from 10 to 200 particles, demonstrating stable and accurate performance even in high particle-density environments.

In the second part of the experiments, the AIFM model was applied to flow field reconstruction in more complex scenarios. The model effectively reconstructed both steady and vortex flow fields, with relative errors remaining within acceptable limits. Even in challenging flow environments, such as vortex-dominated flows, the model demonstrated high accuracy and stability, making it highly suitable for real-world applications requiring reliable flow analysis. Similarly, when applied to the T-junction flow field, the model exhibited strong performance, with errors remaining within a \SI{15}{\percent} margin under typical conditions. In more complex configurations with wider branches or stronger turbulence, the relative errors increased but stayed within acceptable ranges, highlighting the model’s robustness and adaptability to varying flow conditions.

Overall, the AIFM model has shown promising results in particle detection and flow field reconstruction. Its ability to accurately reconstruct particle distributions and flow fields, even in complex environments, demonstrates its broad applicability in real-world scenarios. The model’s adaptability to different configurations and its capacity to handle varying particle densities and flow conditions make it a versatile and effective tool for engineering applications. The insights gained from this study are valuable for optimizing model parameters and configurations, ensuring that the AIFM model remains both efficient and accurate in diverse practical situations.

\section{Conclusion}\label{sec_conclusion}
This work validates the feasibility of the AIFM model for three-dimensional flow measurement and demonstrates its accuracy through numerical experiments. In addition to its high accuracy and adaptability in complex hydrodynamic conditions, the AIFM method offers a significant advantage over traditional flow measurement techniques by enabling remote, non-contact monitoring of river and reservoir flows. This makes it a promising tool for hazardous and large-scale water resource monitoring, where conventional methods face limitations due to safety risks, spatial coverage, or measurement resolution. The main contributions of this work are summarized as follows:

\begin{enumerate}
	\item Development of a realistic model: In this study, we have proposed a model that better reflects real-world scenarios, specifically by incorporating the effects of sidelobe reflections. This addition enhances the accuracy of flow field measurements and improves the overall robustness of the AIFM model in practical applications.
	\item Validation of the AIFM model: Extensive numerical experiments confirm the robustness and accuracy of the AIFM model in particle detection and flow velocity measurement. The results demonstrate that the method maintains high precision across diverse configurations, including variations in receiver layout, particle density, and flow field complexity.
	\item Performance in complex flow scenarios: The AIFM model is successfully applied to three-dimensional T-junction flow, highlighting its adaptability and accuracy under challenging conditions. The method achieves stable performance with low measurement errors across a broad range of geometrical and flow parameters.
\end{enumerate}

Building on the promising numerical results, future research will focus on extending the AIFM model to real-world datasets. This includes integrating field measurements from natural river networks and engineered water systems to validate the performance of the method under practical conditions. Challenges such as noise in measurement data, non-ideal particle distributions, and environmental factors will be addressed. Additionally, further optimization of the AIFM model to enhance its computational efficiency and scalability will be explored, aiming to expand its application to larger and more complex flow systems.


\bibliographystyle{plain}
\bibliography{references}

\end{document}